\documentclass[11pt,fleqn]{article}
\usepackage{cite}
\usepackage{amsmath, amsthm, amssymb}
\usepackage{bm}
\usepackage{float,afterpage}
\usepackage{leftidx}
\usepackage{srcltx}
\usepackage{graphicx}

\input{tcilatex}
\begin{document}

\title{\textbf{Kernel Oriented Generator Distribution}}
\author{{\ A. Bekker$^{1}$\thanks{%
Corresponding Author. Email: andriette.bekker@up.ac.za} and M. Arashi$^{1,2}$} \\
%EndAName
\textit{$^1$Department of Statistics, Faculty of Natural and
Agricultural
Sciences, }\\
\vspace{.5cm} \textit{University of Pretoria, Pretoria, 0002,
South Africa }\\
\textit{$^2$Department of Statistics, School of Mathematical
sciences,
Shahrood University,}\\
\vspace{.5cm} \textit{Shahrood, Iran }}
\date{}
\maketitle

\vspace{0.5cm}

\begin{quotation}
\noindent \textit{Abstract:} Matrix variate beta (MVB) distributions are
used in different fields of hypothesis testing, multivariate correlation
analysis, zero regression, canonical correlation analysis and etc. In this
approach a unified methodology is proposed to generate matrix variate
distributions by combining the kernel of MVB distributions of different
types with an unknown Borel measurable function of trace operator over
matrix space, called generator component. The latter component is a
principal element of these newly defined generator type matrix variate
distributions. The matrix variate Kummer beta distribution is amongst others
a special case. Several statistical properties of this newly defined family
of distributions are derived. In the conclusion other extensions and
developments are discussed.


\vspace{9pt} \noindent \textit{AMS Classification:} Primary:
62H10, Secondary: 62H15
\end{quotation}

%%%%%%%%%%%%%%%%%%%%%%%%%%%%%%%%%%%%%%%%%%%%%%% Greek %%%%%%%%%%%%%%%%%%%%%%%%%%%%
%\newcommand{\bXi}{\boldsymbol{\Xi}}
%%%%%%%%%%%%%%%%%%%%%%%%%%%%%%%%%%%%%%%%%%%%%%%%%%%%%%%%%%%% Alphabets %%%%%%%%%%%%%%%%%%%%%%%%
%\newcommand{\bA}{\boldsymbol{A}}
%\newcommand{\bf}{\boldsymbol{f}}
\renewcommand{\bm}{\boldsymbol{m}}
%%%%%%%%%%%%%%%%%%%%%%%%%%%%%%%%%%%%%%%%%%%%%%%%%%  Math %%%%%%%%%%%%%%%%%%%%%%%%%%%%%
\renewcommand{\H}{\mathcal{H}} \renewcommand{\L}{\mathcal{L}} %
\renewcommand{\O}{\mathcal{O}} \renewcommand{\P}{\mathcal{P}} %
\renewcommand{\S}{\mathcal{S}}
%%%%%%%%%%%%%%%%%%%%%%%%%%%%%%%%%%%%%%%%%%%%%%%%% Theorems %%%%%%%%%%%%%%%%%

\def\blackbox{\ \rule{0.5em}{0.5em}}
\def\tr{\mathop{\rm tr}\nolimits}
\def\Re{\mathop{\rm Re}\nolimits}
\def\Retr{\mathop{\rm Retr}\nolimits}
\def\diag{\mathop{\rm diag}\nolimits}
\def\tr{\mathop{\rm tr}\nolimits}
\def\rank{\mathop{\rm rank}\nolimits}
\def\vec{\mathop{\rm vec}\nolimits}
\def\vecp{\mathop{\rm vecp}\nolimits}
\def\vol{\mathop{\rm Vol}\nolimits}
\def\etr{\mathop{\rm etr}\nolimits}
\def\Rel{\mathop{\rm Re}\nolimits}
\def\var{\mathop{\rm Var}\nolimits}
\def\cov{\mathop{\rm Cov}\nolimits}
\def\corr{\mathop{\rm corr}\nolimits}
\def\func #1{\mathop{\rm #1}\nolimits}%

%%%%%%%%%%%%%%%%%%%%%%%%%%%%%%%%%%%%%%%%%%%%%%% Greek %%%%%%%%%%%%%%%%%%%%%%%%%%%%
\newcommand{\balpha}{\boldsymbol{\alpha}}
\newcommand{\bbeta}{\boldsymbol{\beta}}
\newcommand{\bgamma}{\boldsymbol{\gamma}}
\newcommand{\bdelta}{\boldsymbol{\delta}}
\newcommand{\bepsilon}{\boldsymbol{\epsilon}}
\newcommand{\bvarepsilon}{\boldsymbol{\varepsilon}}
\newcommand{\bzeta}{\boldsymbol{\zeta}}
\newcommand{\bet}{\boldsymbol{\eta}}
\newcommand{\btheta}{\boldsymbol{\theta}}
\newcommand{\biota}{\boldsymbol{\iota}}
\newcommand{\bkappa}{\boldsymbol{\kappa}}
\newcommand{\blambda}{\boldsymbol{\lambda}}
\newcommand{\bmu}{\boldsymbol{\mu}}
\newcommand{\bnu}{\boldsymbol{\nu}}
\newcommand{\bxi}{\boldsymbol{\xi}}
\newcommand{\0}{\boldsymbol{0}}
\newcommand{\1}{\boldsymbol{1}}
\newcommand{\bpi}{\boldsymbol{\pi}}
\newcommand{\bvarpi}{\boldsymbol{\varpi}}
\newcommand{\brho}{\boldsymbol{\rho}}
\newcommand{\bvarrho}{\boldsymbol{\varrho}}
\newcommand{\bsigma}{\boldsymbol{\sigma}}
\newcommand{\bvarsigma}{\boldsymbol{\varsigma}}
\newcommand{\btau}{\boldsymbol{\tau}}
\newcommand{\bupsilon}{\boldsymbol{\upsilon}}
\newcommand{\bphi}{\boldsymbol{\phi}}
\newcommand{\bvarphi}{\boldsymbol{\varphi}}
\newcommand{\bchi}{\boldsymbol{\chi}}
\newcommand{\bpsi}{\boldsymbol{\psi}}
\newcommand{\bomega}{\boldsymbol{\omega}}
\newcommand{\bGamma}{\boldsymbol{\Gamma}}
\newcommand{\bDelta}{\boldsymbol{\Delta}}
\newcommand{\bTheta}{\boldsymbol{\Theta}}
\newcommand{\bLambda}{\boldsymbol{\Lambda}}
\newcommand{\bXi}{\boldsymbol{\Xi}}
\newcommand{\bSigma}{\boldsymbol{\Sigma}}
\newcommand{\bUpsilon}{\boldsymbol{\Upsilon}}
\newcommand{\bPhi}{\boldsymbol{\Phi}}
\newcommand{\bPsi}{\boldsymbol{\Psi}}
\newcommand{\bOmega}{\boldsymbol{\Omega}}
%%%%%%%%%%%%%%%%%%%%%%%%%%%%%%%%%%%%%%%%%%%%%%%%%%%%%%%%%%%% Alphabets %%%%%%%%%%%%%%%%%%%%%%%%
\newcommand{\ba}{\boldsymbol{a}}
\newcommand{\bb}{\boldsymbol{b}}
\newcommand{\bA}{\boldsymbol{A}}
\newcommand{\bB}{\boldsymbol{B}}
\newcommand{\bc}{\boldsymbol{c}}
\newcommand{\bC}{\boldsymbol{C}}
\newcommand{\bd}{\boldsymbol{d}}
\newcommand{\bD}{\boldsymbol{D}}
\newcommand{\be}{\boldsymbol{e}}
\newcommand{\bE}{\boldsymbol{E}}
\newcommand{\bF}{\boldsymbol{F}}
\newcommand{\bg}{\boldsymbol{g}}
\newcommand{\bG}{\boldsymbol{G}}
\newcommand{\bh}{\boldsymbol{h}}
\newcommand{\bH}{\boldsymbol{H}}
\newcommand{\bi}{\boldsymbol{i}}
\newcommand{\bI}{\boldsymbol{I}}
\newcommand{\bj}{\boldsymbol{j}}
\newcommand{\bJ}{\boldsymbol{J}}
\newcommand{\bk}{\boldsymbol{k}}
\newcommand{\bK}{\boldsymbol{K}}
\newcommand{\bl}{\boldsymbol{l}}
\newcommand{\bL}{\boldsymbol{L}}
\renewcommand{\bm}{\boldsymbol{m}}
\newcommand{\bM}{\boldsymbol{M}}
\newcommand{\bn}{\boldsymbol{n}}
\newcommand{\bN}{\boldsymbol{N}}
\newcommand{\bo}{\boldsymbol{o}}
\newcommand{\bO}{\boldsymbol{O}}
\newcommand{\bp}{\boldsymbol{p}}
\newcommand{\bP}{\boldsymbol{P}}
\newcommand{\bq}{\boldsymbol{q}}
\newcommand{\bQ}{\boldsymbol{Q}}
\newcommand{\br}{\boldsymbol{r}}
\newcommand{\bR}{\boldsymbol{R}}
\newcommand{\bs}{\boldsymbol{s}}
\newcommand{\bS}{\boldsymbol{S}}
\newcommand{\bt}{\boldsymbol{t}}
\newcommand{\bT}{\boldsymbol{T}}
\newcommand{\bu}{\boldsymbol{u}}
\newcommand{\bU}{\boldsymbol{U}}
\newcommand{\bv}{\boldsymbol{v}}
\newcommand{\bV}{\boldsymbol{V}}
\newcommand{\bw}{\boldsymbol{w}}
\newcommand{\bW}{\boldsymbol{W}}
\newcommand{\bx}{\boldsymbol{x}}
\newcommand{\bX}{\boldsymbol{X}}
\newcommand{\by}{\boldsymbol{y}}
\newcommand{\bY}{\boldsymbol{Y}}
\newcommand{\bz}{\boldsymbol{z}}
\newcommand{\bZ}{\boldsymbol{Z}}
%%%%%%%%%%%%%%%%%%%%%%%%%%%%%%%%%%%%%%%%%%%%%%%%%%  Math %%%%%%%%%%%%%%%%%%%%%%%%%%%%%
\newcommand{\A}{\mathcal{A}}
\newcommand{\B}{\mathcal{B}}
\newcommand{\C}{\mathcal{C}}
\newcommand{\D}{\mathcal{D}}
\newcommand{\E}{\mathcal{E}}
\newcommand{\F}{\mathcal{F}}
\newcommand{\G}{\mathcal{G}}
\renewcommand{\H}{\mathcal{H}}
\newcommand{\I}{\mathcal{I}}
\newcommand{\J}{\mathcal{J}}
\newcommand{\K}{\mathcal{K}}
\renewcommand{\L}{\mathcal{L}}
\newcommand{\M}{\mathcal{M}}
\newcommand{\N}{\mathcal{N}}
\renewcommand{\O}{\mathcal{O}}
\renewcommand{\P}{\mathcal{P}}
\newcommand{\Q}{\mathcal{Q}}
\newcommand{\R}{\mathcal{R}}
\renewcommand{\S}{\mathcal{S}}
\newcommand{\T}{\mathcal{T}}
\newcommand{\U}{\mathcal{U}}
\newcommand{\V}{\mathcal{V}}
\newcommand{\W}{\mathcal{W}}
\newcommand{\X}{\mathcal{X}}
\newcommand{\Y}{\mathcal{Y}}
\newcommand{\Z}{\mathcal{Z}}
%%%%%%%%%%%%%%%%%%%%%%%%%%%%%%%%%%%%%%%%%%%%%%%%% Theorems %%%%%%%%%%%%%%%%%
\def\blackbox{\ \rule{0.5em}{0.5em}}
\def\tr{\mathop{\rm tr}\nolimits}
\def\Re{\mathop{\rm Re}\nolimits}
\def\Retr{\mathop{\rm Retr}\nolimits}
\def\diag{\mathop{\rm diag}\nolimits}
\def\tr{\mathop{\rm tr}\nolimits}
\def\rank{\mathop{\rm rank}\nolimits}
\def\vec{\mathop{\rm vec}\nolimits}
\def\vecp{\mathop{\rm vecp}\nolimits}
\def\vol{\mathop{\rm Vol}\nolimits}
\def\etr{\mathop{\rm etr}\nolimits}
\def\Rel{\mathop{\rm Re}\nolimits}
\def\var{\mathop{\rm Var}\nolimits}
\def\cov{\mathop{\rm Cov}\nolimits}
\def\corr{\mathop{\rm corr}\nolimits}
\def\func #1{\mathop{\rm #1}\nolimits}%

%%%%%%%%%%%%%%%%%%%%%%%%%%%%%%%%%%%%%%%%%%%%%%%%%%%%%%%%%%%%%%%%%%%%%%%%%%%%%

%\newtheorem{theorem}{Theorem}[section] \newtheorem{remark}{Remark}[section] %
%\newtheorem{definition}{Definition}[section] %
%\newtheorem{corollary}{Corollary}[theorem] \newtheorem{lemma}{Lemma}[section]
%\newenvironment{proof}[1][Proof]{\textbf{#1:} }{\ \rule{0.5em}{0.5em}\medskip\par}
\renewcommand{\thetheorem}{\arabic{section}.\arabic{theorem}} %
\renewcommand{\thedefinition}{\arabic{section}.\arabic{definition}} %
\renewcommand{\thelemma}{\arabic{section}.\arabic{lemma}} %
\renewcommand{\thesection}{\arabic{section}} \renewcommand{\thesubsection}{%
\arabic{section}.\arabic{subsection}} \renewcommand{\thecorollary}{%
\arabic{section}.\arabic{theorem}.\arabic{corollary}}

\renewcommand{\thefootnote}{\fnsymbol{footnote}} %\input{tcilatex}

\section{Introduction}

It is well-documented that change to the stucture of a known
statistical distribution generates a new mutated distribution
which performs better in some cases. One interesting and
well-known approach is to incorporate the kernel of a statistical
distribution to propose another one. Examples include the works
of, but not restricted to, Jones (2004), Nadarajah and Kotz (2004,
2006), Brown et al. (2002), Pauw et al. (2010), Silva et al.
(2010), Singla et al. (2012) and Ferreira et al. (2012). The
weighted distribution is nothing but a mathematical construct to
the statistical distribution where there is usually an underlying
`chance mechanism' associated with the population of interest
(e.g. Nanda and Jain, 1999; Navarro et al., 2006; Kwam, 2008 and
Sunoj and Linu, 2012).

In this paper the authors propose a new kernel-generator
definition that is a composition of a kernel of a statistical
distribution combined with a Borel measurable function of trace
operator over matrix space. The kernel oriented generator
approach, from matrix variate viewpoint, is defined as follows:
\begin{definition}\label{definnition kernel oriented}
The random symmetric matrix $\bX$ has kernel oriented distribution
if it can be represented as
\begin{eqnarray*}
f(\bX)=C_o f^*(\bX)h(\tr\bPhi\bX),\quad \bX>\0,
\end{eqnarray*}
where $f^*(.)$ is the kernel of any statistical distribution,
$h(.)$ is a Borel measurable function which admits Taylor's series
expansion, $\bPhi>\0$ is the canonical parameter and $C_o$ is the
normalizing constant.
\end{definition}We turn the reader's attention to the following:

\begin{enumerate}
\item Note: call $f^{\ast }(.)$ and $h(.)$ as the na\"{\i}ve kernel (NK) and
principal kernel (PK), respectively. \ The latter is also called the
generator.

\item[2.] It would of major task to find the normalizing constant $C_{o}$,
since the PK component can be any function. Recall that an elliptically
contoured distribution (even matrix variate form) is a distribution whose
characteristic function (density if exists) can be presented as a function
of quadratic forms. Thus there is a similarity for the constant $C_{o}$ in
the literature. However, we will address the solution to obtain $C_{o}$ by
applying the Taylor's series expansion under some mild regularity conditions.

\item[3.] It is thus possible to extend each statistical distribution by
taking its NK component and compose it with an extra PK element, which gives
an infinity class of distributions. The latter element has many statistical
features where the shape of $f(.)$ is the important one. This approach is
not restricted to matrix variate distributions only, however some univariate
examples are also considered here.
\end{enumerate}

Using Definition \ref{definnition kernel oriented} in this paper,
we focus on well-known matrix variate beta kernels. The resulting
new distributions will be referred to as matrix variate beta
kernel oriented generator distributions or matrix variate beta
type 1/2/3 generator distribution (MBG1/2/3) for short.

We organize the paper as follows: In section 2 the definitions of
\ the matrix variate beta generator distributions of type I, II
and III are given. Section 3 is devoted to some important
statistical properties of these new distributions, followed by a
discussion section. The expressions are given in terms of zonal
polynomials, homogeneous invariant polynomials with two or more
matrix arguments, Meijer's G function. The reader is referred to
the papers of (Chikuse, 1980; Davis, 1979, 1980 and James, 1961,
1964).

\section{Matrix Variate Beta Generator Distributions}

The well-known matrix variate beta distributions (Olkin and Rubin,
1962), used in different fields of hypothesis testing,
multivariate correlation analysis, zero regression, have been
extended by several authors. The matrix variate beta type 3
distribution has been defined, and some of its properties have
been studied by Gupta and Nagar (2000b, 2009). More recently Nagar
et al. (2013), by using extended matrix variate beta function,
generalized the well-known matrix beta type 1 distribution. Gupta
and Nagar (2006) extended the work of Nadarajah and Kotz (2006) by
defining matrix variate hypergeometric beta distribution. Ehlers
(2011) proposed the matrix variate beta type 5 distribution
motivating from generalized hypothesis testing in multivariate
setup (see also Bekker et al., 2012).\bigskip

Let $\boldsymbol{X}$ be a random symmetric matrix of dimension $m$ and $%
\mathop{\rm Re}\nolimits(a),\mathop{\rm Re}\nolimits(b)>(m-1)/2$. According
to Definition 1.1, in this section we define

\begin{enumerate}
\item[(i)] matrix variate beta type 1 generator distribution (MBG1), by
taking the NK to be
\begin{equation*}
\det(\boldsymbol{X})^{a-\frac{1}{2}(m+1)}\det(\boldsymbol{I}-\boldsymbol{X}%
)^{b-\frac{1}{2}(m+1)}
\end{equation*}

\item[(ii)] matrix variate beta type 2 generator distribution (MBG2), by
taking the NK to be
\begin{equation*}
\det(\boldsymbol{X})^{a-\frac{1}{2}(m+1)}\det(\boldsymbol{I}+\boldsymbol{X}%
)^{-(a+b)}
\end{equation*}

\item[(iii)] matrix variate beta type 3 generator distribution (MBG3), by
taking the NK to be
\begin{equation*}
\det(\boldsymbol{X})^{a-\frac{1}{2}(m+1)}\det(\boldsymbol{I}-\boldsymbol{X}%
)^{b-\frac{1}{2}(m+1)}\det(\boldsymbol{I}+\boldsymbol{X})^{-(a+b)}
\end{equation*}
\end{enumerate}

Further we consider some special cases.

\begin{definition}\label{MBG}
The random symmetric matrix $\bX$ of dimension $m$ is said to have
\begin{enumerate}
\item[(i)] MBG1 distribution with parameters
$a$, $b$ and $\bPhi$ and shape generator $h$, denoted by $\bX\sim
MBG1_m(a,b,h,\bPhi)$, if it has the following density function
\begin{eqnarray*}
f(\bX)=\zeta_{a,b}^{(1)}\det(\bX)^{a-\frac{1}{2}(m+1)}\det(\bI_m-\bX)^{b-\frac{1}{2}(m+1)}h(\tr(\bPhi\bX)),
\quad \bX\in\I_m,
\end{eqnarray*}
\item[(ii)] MBG2 distribution with parameters
$a$, $b$ and $\bPhi$ and shape generator $h$, denoted by $\bX\sim
MBG2_m(a,b,h,\bPhi)$, if it has the following density function
\begin{eqnarray*}
f(\bX)=\zeta_{a,b}^{(2)}\det(\bX)^{a-\frac{1}{2}(m+1)}\det(\bI_m+\bX)^{-(a+b)}h(\tr(\bPhi\bX)),
\quad \bX\in\S_m,
\end{eqnarray*}
\item[(iii)] MBG3 distribution with parameters
$a$, $b$ and $\bPhi$ and shape generator $h$, denoted by
$\bX\sim MBG3_m(a,b,h,\bPhi)$, if it has the following density
function
\begin{eqnarray*}
f(\bX)=\zeta_{a,b}^{(3)}\det(\bX)^{a-\frac{1}{2}(m+1)}\det(\bI_m-\bX)^{b-\frac{1}{2}(m+1)}\det(\bI_m+\bX)^{-(a+b)}h(\tr(\bPhi\bX)),
\quad \bX\in\I_m,
\end{eqnarray*}
\end{enumerate}
where, $\S_m$ is the space of all positive definite matrices of order $m$, $\I_m$ is the space of all square matrices of order $m$ such that $\bI_m-\bX\in\S_m$ iff $\bX\in\S_m$, and $\Re(a)>(m-1)/2$, $\Re(b)>(m-1)/2$, $\bPhi$ is a
symmetric complex matrix, $h(.)$ is a Borel measurable function
that admits a Taylor's series expansion and $\zeta_{a,b}^{(j)}$,
$j=1,2,3$ are the normalizing constants.
\end{definition}%
\begin{remark}
To find the normalizing constants in Definition \ref{MBG}, first
we use the Taylor's series expansion to get
\begin{equation}
h(\tr\bPhi\bX)=\sum_{t=0}^\infty
\frac{h^{(t)}(0)}{t!}\;\tr(\bPhi\bX)^t=\sum_{t=0}^\infty
\frac{h^{(t)}(0)}{t!}\sum_\tau C_\tau(\bPhi\bX),
\end{equation}
where $C_\tau(.)$ is the zonal polynomial, and we used ordered partitions in use of
the zonal polynomials. Then
$\zeta_{a,b}^{(j)}$, $j=1,2,3$ can be obtained after some matrix
algebra as:
\begin{eqnarray*}
\left(\zeta_{a,b}^{(1)}\right)^{-1} &=& \sum_{t=0}^\infty
                          \frac{h^{(t)}(0)}{t!}\sum_\tau\int_{\I_m}
\det(\bX)^{a-\frac{1}{2}(m+1)}\det(\bI_m-\bX)^{b-\frac{1}{2}(m+1)}C_\tau(\bPhi\bX)\textnormal{d}\bX\cr
                          &=& \sum_{t=0}^\infty
                          \frac{h^{(t)}(0)}{t!}\sum_\tau\frac{\Gamma_m(a,\tau)\Gamma_m(b)}{\Gamma_m(a+b,\tau)}C_\tau(\bPhi),\quad\mbox{Theorem 7.2.10\;of\;Muirhead\;(2005)}\cr
\left(\zeta_{a,b}^{(2)}\right)^{-1} &=& \sum_{t=0}^\infty
                          \frac{h^{(t)}(0)}{t!}\sum_\tau\int_{\S_m}
\det(\bX)^{a-\frac{1}{2}(m+1)}\det(\bI_m+\bX)^{-(a+b)}C_\tau(\bPhi\bX)\textnormal{d}\bX\cr
                          &=& \sum_{t=0}^\infty
                          \frac{h^{(t)}(0)}{t!}\sum_\tau\frac{\Gamma_m(a,\tau)\Gamma_m(b,-\tau)}{\Gamma_m(a+b)}C_\tau(\bPhi),\quad\mbox{Lemma
                          5\;of\;Khatri\;(1966)}\cr
\left(\zeta_{a,b}^{(3)}\right)^{-1} &=& \sum_{t=0}^\infty
                          \frac{h^{(t)}(0)}{t!}\sum_\tau\int_{\I_m}
\det(\bX)^{a-\frac{1}{2}(m+1)}\det(\bI_m-\bX)^{b-\frac{1}{2}(m+1)}\det(\bI_m+\bX)^{-(a+b)}C_\tau(\bPhi\bX)\textnormal{d}\bX\cr
                          &=& \sum_{t=0}^\infty
                          \frac{h^{(t)}(0)}{t!}\sum_\tau\sum_{k=0}^\infty\sum_\kappa\frac{(-1)^k(a+b)_\kappa}{k!}\cr
                          &&\int_{\I_m}
\det(\bX)^{a-\frac{1}{2}(m+1)}\det(\bI_m-\bX)^{b-\frac{1}{2}(m+1)}C_\kappa(\bX)C_\tau(\bPhi\bX)\textnormal{d}\bX\quad\cr
                          &=& \sum_{t=0}^\infty
                          \frac{h^{(t)}(0)}{t!}\sum_\tau\sum_{k=0}^\infty\sum_\kappa\frac{(-1)^k(a+b)_\kappa}{k!}\sum_{\phi\in\kappa\cdot\tau}\theta_\phi^{\kappa,\tau}
\cr
                          &&\int_{\I_m}
\det(\bX)^{a-\frac{1}{2}(m+1)}\det(\bI_m-\bX)^{b-\frac{1}{2}(m+1)}C_\phi^{\kappa,\tau}(\bX,\bPhi\bX)\textnormal{d}\bX,
\quad\mbox{Eq. (2.8)\; of\; Davis\;(1979)}\cr
                          &=& \sum_{t=0}^\infty
                          \frac{h^{(t)}(0)}{t!}\sum_\tau\sum_{k=0}^\infty\sum_\kappa\frac{(-1)^k(a+b)_\kappa}{k!}\sum_{\phi\in\kappa\cdot\tau}\theta_\phi^{\kappa,\tau}
                          \cr
                          &&\frac{\Gamma_m(b)\Gamma_m(a,\phi)}{\Gamma_m(a+b,\phi)}C_\phi^{\kappa,\tau}(\bI_m,\bPhi),\quad\mbox{Eq.
(3.28)\; of\; Chikuse\;(1980)}\cr
                          &=&\sum_{\tau,\kappa,\phi}\frac{(-1)^k(a+b)_\kappa h^{(t)}(0)}{t!k!}
                          \frac{\Gamma_m(b)\Gamma_m(a,\phi)}{\Gamma_m(a+b,\phi)}\frac{\left(\theta_\phi^{\kappa,\tau}\right)^2C_\phi(\bI_m)}{C_\tau(\bI_m)}C_\tau(\bPhi),\cr
                          &&\qquad\qquad\qquad\qquad\qquad\qquad\qquad\qquad\qquad\qquad\qquad\qquad\mbox{Eq.
(2.2)\; of\; Davis\;(1979)},
\end{eqnarray*}
where
$\theta_\phi^{\kappa,\tau}=C_\phi^{\kappa,\tau}(\bI_m,\bI_m)/C_\phi(\bI_m)$,
$
\sum_{\tau,\kappa,\phi}\equiv\sum_{t=0}^\infty\sum_{k=0}^\infty\sum_\tau\sum_\kappa\sum_{\phi\in\kappa\cdot\tau}
$, $\Gamma_p(.)$ represents the multivariate gamma function, and
$\Gamma_p(.,\kappa)$ the generalized gamma function of weight
$\kappa$. (See Gupta and Nagar, 2000a)
\end{remark}

\section{Characteristics}

In this section we provide some important statistical properties for three
different types of matrix variate beta generator distributions.

The following result is straightforward.
\begin{theorem}\label{expectation}
Let $\bX_i\sim MBGi_m(a,b,h,\bPhi)$, $i=1,2,3$. Then it follows
that
\begin{eqnarray*}
E\left(\det(\bX_1)^r\right) &=& \zeta_{a,b}^{(1)}\sum_{t=0}^\infty
                          \frac{h^{(t)}(0)}{t!}\sum_\tau\frac{\Gamma_m(a+r,\tau)\Gamma_m(b)}{\Gamma_m(a+r+b,\tau)}C_\tau(\bPhi),\cr
E\left(\det(\bX_2)^r\right) &=& \zeta_{a,b}^{(2)}\sum_{t=0}^\infty
                          \frac{h^{(t)}(0)}{t!}\sum_\tau\frac{\Gamma_m(a+r,\tau)\Gamma_m(b,-\tau)}{\Gamma_m(a+r+b)}C_\tau(\bPhi),\cr
E\left(\det(\bX_3)^r\right) &=&
\zeta_{a,b}^{(3)}\sum_{\tau,\kappa,\phi}\frac{(-1)^k(a+b)_\kappa
h^{(t)}(0)}{t!k!}
\frac{\Gamma_m(a+r,\phi)\Gamma_m(b)}{\Gamma_m(a+b+r,\phi)}\frac{\left(\theta_\phi^{\kappa,\tau}\right)^2C_\phi(\bI_m)}{C_\tau(\bI_m)}C_\tau(\bPhi)
\end{eqnarray*}
\end{theorem}

In the following theorem, we give the moment generating function (MGF) for
each type of MBG distribution.

\begin{theorem}\label{MGF}
Denote the MGF of $\bX_i\sim MBGi_m(a,b,h,\bPhi)$, $i=1,2,3$ by
$\M_i$. Then we have
\begin{eqnarray*}
\M_1(\bT)&=&\zeta_{a,b}^{(1)}\sum_{\tau,\kappa,\phi}\frac{h^{(t)}(0)}{t!k!}\theta_\phi^{\kappa,\tau}
\frac{\Gamma_m(a,\phi)\Gamma_m(b)
}{\Gamma_m(a+b,\phi)}C_\phi^{\kappa,\tau}(\bT,\bPhi),\cr
\M_2(\bT)&=&\zeta_{a,b}^{(2)}\det(\bT)^{-a}\sum_{\tau,\kappa,\phi}\frac{h^{(t)}(0)}{t!}\frac{(-1)^{-(am+t)}(a+b)_\kappa}{k!}
\theta_\phi^{\kappa,\tau}
\Gamma_m(a,\phi)C_\phi(\bT^{-1},\bT^{-1}\bPhi),\cr
\M_3(\bT)&=&\zeta_{a,b}^{(3)}
\sum_{\tau,\kappa,\lambda,\phi}\frac{h^{(t)}(0)}{t!k!}\frac{(-1)^l(a+b)_\lambda}{l!}\theta_\phi^{\kappa,\tau,\lambda}
\frac{\Gamma_m(a,\phi)\Gamma_m(b)}{\Gamma_m(a+b,\phi)}C_\phi^{\kappa,\tau,\lambda}(\bT,\bPhi,\bI),
\end{eqnarray*}
where
$\sum_{\tau,\kappa,\lambda,\phi}\equiv\sum_{t=0}^\infty\sum_\tau\sum_{k=0}^\infty\sum_\kappa\sum_{l=0}^\infty\sum_\lambda\sum_{\phi\in\kappa\cdot\tau\cdot\lambda}$
and
$\theta_\phi^{\kappa,\tau,\lambda}=C_\phi^{\kappa,\tau,\lambda}(\bI_m,\bI_m,\bI_m)/C_\phi(\bI_m)$.
\end{theorem}
\noindent\textbf{Proof:} The proof of $\mathcal{M}_1$ is straightforward.
Here we provide the proofs of $\mathcal{M}_2$ \& $\mathcal{M}_3$.

Using Taylor's series expansion for $h$, and Eq. (3.10) of Chikuse (1980) we
get
\begin{eqnarray*}
\mathcal{M}_2(\boldsymbol{T})&=&\zeta_{a,b}^{(2)}\int_{\S _m}\det(%
\boldsymbol{X})^{a-\frac{1}{2}(m+1)}\det(\boldsymbol{I}_m+\boldsymbol{X}%
)^{-(a+b)}\cr &&h(\mathop{\rm tr}\nolimits(\boldsymbol{\Phi}\boldsymbol{X}))%
\mathop{\rm etr}\nolimits(\boldsymbol{T}\boldsymbol{X}) \mathnormal{d}%
\boldsymbol{X} \cr &=&\zeta_{a,b}^{(2)}\sum_{t=0}^\infty\sum_\tau
\sum_{k=0}^\infty\sum_\kappa\frac{h^{(t)}(0)}{t!}\frac{(-1)^k(a+b)_\kappa}{k!%
}\cr && \int_{\S _m}\det(\boldsymbol{X})^{a-\frac{1}{2}(m+1)}C_\tau(%
\boldsymbol{\Phi}\boldsymbol{X})C_\kappa(\boldsymbol{X})\etr(\boldsymbol{T}%
\boldsymbol{X})\mathnormal{d}\boldsymbol{X}\cr &=&\zeta_{a,b}^{(2)}%
\sum_{t=0}^\infty\sum_\tau \sum_{k=0}^\infty\sum_\kappa\frac{h^{(t)}(0)}{t!}%
\frac{(-1)^k(a+b)_\kappa}{k!} \sum_{\phi\in\kappa\cdot\tau}\theta_\phi^{%
\kappa,\tau}\cr && \int_{\S _m}\det(\boldsymbol{X})^{a-\frac{1}{2}%
(m+1)}C_\phi^{\kappa,\tau}(\boldsymbol{X},\boldsymbol{\Phi}\boldsymbol{X}%
)\etr(\boldsymbol{T}\boldsymbol{X})\mathnormal{d}\boldsymbol{X}.
\notag
\end{eqnarray*}
Using Eq. (3.21) of Chikuse (1980) we obtain
\begin{eqnarray*}
\mathcal{M}_2(\boldsymbol{T})&=&\zeta_{a,b}^{(2)}\sum_{t=0}^\infty\sum_\tau
\sum_{k=0}^\infty\sum_\kappa\frac{h^{(t)}(0)}{t!}\frac{(-1)^k(a+b)_\kappa}{k!%
} \sum_{\phi\in\kappa\cdot\tau}\theta_\phi^{\kappa,\tau}\cr &&
\Gamma_m(a,\phi)(-1)^{-(am+t+k)}\det(\boldsymbol{T})^{-a}C_\phi(\boldsymbol{T}%
^{-1},\boldsymbol{T}^{-1}\boldsymbol{\Phi}).
\end{eqnarray*}
For the MGF of MBG3, by making use of Taylor's series expansion for $h$, and
Eq. (3.10) of Chikuse (1980) we get %\begin{eqnarray}\label{Taylor}
%h(\tr(\bPhi\bX))=\sum_{t=0}^\infty
%\frac{h^{(t)}(0)}{t!}[\tr(\bPhi\bX)]^t=\sum_{t=0}^\infty
%\frac{h^{(t)}(0)}{t!}\sum_\tau C_\tau(\bPhi\bX).
%\end{eqnarray}

\begin{eqnarray*}
\mathcal{M}_3(\boldsymbol{T})&=&\zeta_{a,b}^{(3)}\int_{\mathcal{I}_m}\det(%
\boldsymbol{X})^{a-\frac{1}{2}(m+1)}\det(\boldsymbol{I}_m-\boldsymbol{X})^{b-%
\frac{1}{2}(m+1)}\det(\boldsymbol{I}_m+\boldsymbol{X})^{-(a+b)}\cr &&h(%
\mathop{\rm tr}\nolimits(\boldsymbol{\Phi}\boldsymbol{X}))\mathop{\rm etr}%
\nolimits(\boldsymbol{T}\boldsymbol{X}) \mathnormal{d}\boldsymbol{X} \cr %
&=&\zeta_{a,b}^{(3)}\sum_{t=0}^\infty\sum_\tau \sum_{k=0}^\infty\sum_\kappa%
\frac{h^{(t)}(0)}{t!k!}\sum_{l=0}^\infty\sum_\lambda\frac{(-1)^l(a+b)_\lambda%
}{l!}\cr && \int_{\mathcal{I}_m}\det(\boldsymbol{X})^{a-\frac{1}{2}%
(m+1)}\det(\boldsymbol{I}_m-\boldsymbol{X})^{b-\frac{1}{2}(m+1)}C_\tau(%
\boldsymbol{\Phi}\boldsymbol{X})C_\kappa(\boldsymbol{T}\boldsymbol{X}%
)C_\lambda(\boldsymbol{X})\mathnormal{d}\boldsymbol{X} \cr %
&=&\zeta_{a,b}^{(3)}\sum_{t=0}^\infty\sum_\tau \sum_{k=0}^\infty\sum_\kappa%
\frac{h^{(t)}(0)}{t!k!}\sum_{l=0}^\infty\sum_\lambda\frac{(-1)^l(a+b)_\lambda%
}{l!} \sum_{\phi\in\kappa\cdot\tau\cdot\lambda}\theta_\phi^{\kappa,\tau,%
\lambda}\cr && \int_{\mathcal{I}_m}\det(\boldsymbol{X})^{a-\frac{1}{2}%
(m+1)}\det(\boldsymbol{I}_m-\boldsymbol{X})^{b-\frac{1}{2}%
(m+1)}C_\phi^{\kappa,\tau,\lambda}(\boldsymbol{T}\boldsymbol{X},\boldsymbol{%
\Phi}\boldsymbol{X},\boldsymbol{X})\mathnormal{d}\boldsymbol{X}.  \notag
\end{eqnarray*}
Finally applying Eq. (3.28) of Chikuse (1980), yields
\begin{eqnarray*}
\mathcal{M}_3(\boldsymbol{T})&=& \zeta_{a,b}^{(3)}\sum_{t=0}^\infty\sum_\tau
\sum_{k=0}^\infty\sum_\kappa\frac{h^{(t)}(0)}{t!k!}\sum_{l=0}^\infty\sum_%
\lambda\frac{(-1)^l(a+b)_\lambda}{l!}\cr &&
\sum_{\phi\in\kappa\cdot\tau\cdot\lambda}\theta_\phi^{\kappa,\tau,\lambda}
\frac{\Gamma_m(a,\phi)\Gamma_m(b)}{\Gamma_m(a+b,\phi)}C_\phi^{\kappa,\tau,%
\lambda}(\boldsymbol{T},\boldsymbol{\Phi},\boldsymbol{I}).
\end{eqnarray*}
and the proof is complete.\hfill$\blacksquare$

In the following result, we give the exact expressions for the cumulative
distribution function (CDF) of MBG1/2/3 distribution.
\begin{theorem}\label{CDF}
Denote the CDF of $\bX_i\sim MBGi_m(a,b,h,\bPhi)$, $i=1,2,3$ by
$F_i$. Then we have
\begin{eqnarray*}
F_1(\bY)&=&\zeta_{a,b}^{(1)}\det(\bY)^{a}\sum_{\tau,\kappa,\phi}\frac{\left(-b+\frac{1}{2}(m+1)\right)_\kappa}{k!}
\frac{h^{(t)}(0)}{t!}
\frac{\Gamma_m(a,\phi)\Gamma\left(\frac{m+1}{2}\right)}{\Gamma\left(a+\frac{m+1}{2},\phi\right)}\cr
&&\times
\theta_\phi^{\kappa,\tau}C_\phi^{\kappa,\tau}(\bY,\bY^{\frac{1}{2}}\bPhi\bY^{\frac{1}{2}}),\cr
F_2(\bY)&=&\zeta_{a,b}^{(2)}\det(\bY)^{a}\sum_{\tau,\kappa,\phi}\frac{(-1)^k\left(a+b\right)_\kappa}{k!}
\frac{h^{(t)}(0)}{t!}
\frac{\Gamma_m(a,\phi)\Gamma\left(\frac{m+1}{2}\right)}{\Gamma\left(a+\frac{m+1}{2},\phi\right)}\cr
&&\times
\theta_\phi^{\kappa,\tau}C_\phi^{\kappa,\tau}(\bY,\bY^{\frac{1}{2}}\bPhi\bY^{\frac{1}{2}}),\cr
F_3(\bY)&=&\zeta_{a,b}^{(3)}\det(\bY)^{a}\sum_{\kappa,\tau,\lambda,\phi}\frac{\left(-b+\frac{1}{2}(m+1)\right)_\kappa}{k!}
\frac{h^{(t)}(0)}{t!}\frac{(-1)^l(a+b)_\lambda}{l!}\frac{\Gamma_m\left(\frac{m+1}{2}\right)\Gamma_m(a,\phi)}{\Gamma_m\left(a+\frac{m+1}{2},\phi\right)}\cr
&&\times \theta_\phi^{\kappa,\tau,\lambda}
C_\phi^{\kappa,\tau,\lambda}(\bY,\bY^{\frac{1}{2}}\bPhi\bY^{\frac{1}{2}},\bY)
\end{eqnarray*}
\end{theorem}\noindent \textbf{Proof:} For the CDF of MBG1 distribution, we
have by definition
\begin{equation*}
F_{1}(\boldsymbol{Y})=\zeta _{a,b}^{(1)}\int_{0<\boldsymbol{X}<\boldsymbol{Y}%
<\boldsymbol{I}_{m}}\det (\boldsymbol{X})^{a-\frac{1}{2}(m+1)}\det (%
\boldsymbol{I}_{m}-\boldsymbol{X})^{b-\frac{1}{2}(m+1)}h(\mathop{\rm tr}%
\nolimits(\boldsymbol{\Phi }\boldsymbol{X}))\mathnormal{d}\boldsymbol{X}.
\end{equation*}%
Making the transformation $\boldsymbol{G}=\boldsymbol{Y}^{-\frac{1}{2}}%
\boldsymbol{X}\boldsymbol{Y}^{-\frac{1}{2}}$ with the Jacobian $J(%
\boldsymbol{X}\rightarrow \boldsymbol{G})=\det (\boldsymbol{Y})^{\frac{1}{2}%
(m+1)}$ we get
\begin{eqnarray*}
F_{1}(\boldsymbol{Y})&=&\zeta _{a,b}^{(1)}\det (\boldsymbol{Y})^{a}\int_{%
\mathcal{I}_{m}}\det (\boldsymbol{G})^{a-\frac{1}{2}(m+1)}\det (\boldsymbol{I%
}_{p}-\boldsymbol{Y}^{\frac{1}{2}}\boldsymbol{G}\boldsymbol{Y}^{\frac{1}{2}%
})^{b-\frac{1}{2}(m+1)}\cr
&&\times h\left( \mathop{\rm tr}\nolimits(%
\boldsymbol{\Phi }\boldsymbol{Y}^{\frac{1}{2}}\boldsymbol{G}\boldsymbol{Y}^{%
\frac{1}{2}})\right) \mathnormal{d}\boldsymbol{G}\cr
&=&\zeta _{a,b}^{(1)}\det (%
\boldsymbol{Y})^{a}\sum_{k=0}^{\infty }\sum_{\kappa }\frac{\left(-b+\frac{1%
}{2}(m+1)\right) _{\kappa }}{k!}\int_{\mathcal{I}_{m}}\det (\boldsymbol{G}%
)^{a-\frac{1}{2}(m+1)}\cr
&&\times C_{\kappa }(\boldsymbol{Y}^{\frac{1}{2}}%
\boldsymbol{G}\boldsymbol{Y}^{\frac{1}{2}})h\left( \mathop{\rm tr}\nolimits(%
\boldsymbol{\Phi }\boldsymbol{Y}^{\frac{1}{2}}\boldsymbol{G}\boldsymbol{Y}^{%
\frac{1}{2}})\right) \mathnormal{d}\boldsymbol{G}.
\end{eqnarray*}%
Make use of Taylor's series expansion for $h(.)$ term and
equations (3.10) and (3.32) of Chikuse (1980) to obtain
\begin{eqnarray*}
F_{1}(\boldsymbol{Y})&=&\zeta _{a,b}^{(1)}\det (\boldsymbol{Y}%
)^{a}\sum_{k=0}^{\infty }\sum_{\kappa }\frac{\left( -b+\frac{1}{2}%
(m+1)\right) _{\kappa }}{k!}\sum_{t=0}^{\infty }\frac{h^{(t)}(0)}{t!}%
\sum_{\tau }\cr
&&\times \int_{\mathcal{I}_{m}}\det (\boldsymbol{G})^{a-\frac{1%
}{2}(m+1)}C_{\phi}^{\kappa,\tau}(\boldsymbol{Y}\boldsymbol{G},\boldsymbol{Y}^{%
\frac{1}{2}}\boldsymbol{\Phi }\boldsymbol{Y}^{\frac{1}{2}}\boldsymbol{G})%
\mathnormal{d}\boldsymbol{G}\cr
&=&\zeta _{a,b}^{(1)}\det (\boldsymbol{Y}%
)^{a}\sum_{k=0}^{\infty }\sum_{\kappa }\frac{\left( -b+\frac{1}{2}%
(m+1)\right) _{\kappa }}{k!}\sum_{t=0}^{\infty }\frac{h^{(t)}(0)}{t!}%
\sum_{\tau }\sum_{\phi \in \kappa \cdot \tau }\theta _{\phi }^{\kappa ,\tau }%
\cr
&&\times \frac{\Gamma _{m}(a,\phi )\Gamma_m \left( \frac{1}{2}(m+1)\right) }{%
\Gamma_m \left( a+\frac{1}{2}(m+1),\phi \right) }C_{\phi }^{\kappa ,\tau }(%
\boldsymbol{Y},\boldsymbol{Y}^{\frac{1}{2}}\boldsymbol{\Phi }\boldsymbol{Y}^{%
\frac{1}{2}}).
\end{eqnarray*}%
The CDF of MBG2 distribution can be obtained in the same fashion
as for the CDF of MBG1 distribution. For the CDF of MBG3
distribution, using the same procedure as in the proof of the CDF
of MBG1 distribution, we have
\begin{eqnarray*}
F_{3}(\boldsymbol{Y})&=&\zeta _{a,b}^{(3)}\det (\boldsymbol{Y}%
)^{a}\sum_{k=0}^{\infty }\sum_{\kappa }\frac{\left( -b+\frac{1}{2}%
(m+1)\right) _{\kappa }}{k!}\sum_{t=0}^{\infty }\frac{h^{(t)}(0)}{t!}%
\sum_{\tau }\sum_{l=0}^{\infty }\sum_{\lambda }\frac{(-1)^{l}(a+b)_{\lambda }%
}{l!}\cr
&&\times \int_{\mathcal{I}_{m}}\det (\boldsymbol{G})^{a-\frac{1}{2}%
(m+1)}C_{\kappa }(\boldsymbol{Y}\boldsymbol{G})C_{\tau }(\boldsymbol{Y}^{%
\frac{1}{2}}\boldsymbol{\Phi }\boldsymbol{Y}^{\frac{1}{2}}\boldsymbol{G}%
)C_{\lambda }(\boldsymbol{Y}\boldsymbol{G})\mathnormal{d}\boldsymbol{G}\cr%
&=&\zeta _{a,b}^{(3)}\det (\boldsymbol{Y})^{a}\sum_{k=0}^{\infty
}\sum_{\kappa
}\frac{\left( -b+\frac{1}{2}(m+1)\right) _{\kappa }}{k!}\sum_{t=0}^{\infty }%
\frac{h^{(t)}(0)}{t!}\sum_{\tau }\sum_{l=0}^{\infty }\sum_{\lambda }\frac{%
(-1)^{l}(a+b)_{\lambda }}{l!}\sum_{\phi \in \kappa \cdot \tau \cdot \lambda
}\theta _{\phi }^{\kappa ,\tau ,\lambda }\cr
&&\times \int_{\mathcal{I}%
_{m}}\det (\boldsymbol{G})^{a-\frac{1}{2}(m+1)}C_{\phi }^{\kappa ,\tau
,\lambda }(\boldsymbol{Y}\boldsymbol{G},\boldsymbol{Y}^{\frac{1}{2}}%
\boldsymbol{\Phi }\boldsymbol{Y}^{\frac{1}{2}}\boldsymbol{G},\boldsymbol{Y}%
\boldsymbol{G})\mathnormal{d}\boldsymbol{G},\quad\mbox{Eq.\;(3.10)\;of\;Chikuse\;(1980)}.
\end{eqnarray*}%
Make use of Eq. (3.32) of Chikuse (1980) to get
\begin{eqnarray*}
F_{3}(\boldsymbol{Y})&=&\zeta _{a,b}^{(3)}\det (\boldsymbol{Y}%
)^{a}\sum_{k=0}^{\infty }\sum_{\kappa }\frac{\left( -b+\frac{1}{2}%
(m+1)\right) _{\kappa }}{k!}\sum_{t=0}^{\infty }\frac{h^{(t)}(0)}{t!}%
\sum_{\tau }\sum_{l=0}^{\infty }\sum_{\lambda }\frac{(-1)^{l}(a+b)_{\lambda }%
}{l!}\sum_{\phi \in \kappa \cdot \tau \cdot \lambda }\theta _{\phi
}^{\kappa ,\tau ,\lambda }\cr &&\times \frac{\Gamma _{m}\left(
\frac{m+1}{2}\right)
\Gamma _{m}(a,\phi )}{\Gamma _{m}\left( a+\frac{m+1}{2},\phi \right) }%
C_{\phi }^{\kappa ,\tau ,\lambda }(\boldsymbol{Y},\boldsymbol{Y}^{\frac{1}{2}%
}\boldsymbol{\Phi }\boldsymbol{Y}^{\frac{1}{2}},\boldsymbol{Y})
\end{eqnarray*}%
which completes the proof.\hfill $\blacksquare $\newline
In what follows, we are interested in the distribution of quadratic forms
from MBG distributions. Assume $\boldsymbol{\Psi },%
\boldsymbol{\Omega }\in \S _{m}$ are some known matrix parameters and under
the meaning of partial l\"{o}wner ordering, $\boldsymbol{\Omega }>%
\boldsymbol{\Psi }$. We are interested in the distribution of the
random matrix variate
\begin{equation}
\boldsymbol{Y}_{i}=(\boldsymbol{\Omega }-\boldsymbol{\Psi })^{\frac{1}{2}}%
\boldsymbol{X}_{i}(\boldsymbol{\Omega }-\boldsymbol{\Psi })^{\frac{1}{2}}+%
\boldsymbol{\Psi },  \label{composite}
\end{equation}%
where $\boldsymbol{X}_{i}\sim MBGi_{m}(a,b,h,\boldsymbol{\Phi })$, $i=1,2,3$%
. The distribution of $\boldsymbol{Y}_{i}$ is the MBG
distribution, which is given in the following result.
\begin{theorem}\label{extended MBG} Suppose that $f(\bY_i)$ is
the density function given by $\eqref{composite}$, while
$\bX_i\sim MBGi_m(a,b,h,\bPhi)$, $i=1,2,3$. Then we have
\begin{eqnarray*}
f(\bY_1)&=&\zeta_{a,b}^{(1)}\det(\bY_1-\bPsi)^{a-\frac{1}{2}(m+1)}\det(\bOmega-\bY_1)^{b-\frac{1}{2}(m+1)}\det(\bOmega-\bPsi)^{-(a+b)+\frac{1}{2}(m+1)}
h(\tr\bTheta(\bY_1-\bPsi)),\cr
f(\bY_2)&=&\zeta_{a,b}^{(2)}\det(\bY_2-\bPsi)^{a-\frac{1}{2}(m+1)}\det(\bOmega+\bY_2-2\bPsi)^{-(a+b)}\det(\bOmega-\bPsi)^{b}
h(\tr\bTheta(\bY_2-\bPsi))\cr
f(\bY_3)&=&\zeta_{a,b}^{(3)}\det(\bY_3-\bPsi)^{a-\frac{1}{2}(m+1)}\det(\bOmega-\bY_3)^{b-\frac{1}{2}(m+1)}\det(\bOmega+\bY_3-2\bPsi)^{-(a+b)}\cr
&&\det(\bOmega-\bPsi)^{\frac{1}{2}(m+1)}
h(\tr\bTheta(\bY_3-\bPsi)),
\end{eqnarray*}
where $\bPsi<\bY_i<\bOmega$, $i=1,2,3$ and
$\bTheta=(\bOmega-\bPsi)^{-\frac{1}{2}}\bPhi(\bOmega-\bPsi)^{-\frac{1}{2}}$.
\end{theorem}\noindent \textbf{Proof:} From Definition \ref{MBG} and the
fact that the Jacobian of transformation is $J(\boldsymbol{X}\rightarrow
\boldsymbol{Y})=\det (\boldsymbol{\Omega }-\boldsymbol{\Psi })^{-\frac{1}{2}%
(m+1)}$, the result follows.\hfill $\blacksquare $.
\begin{remark}
Suppose that $\bX_i\sim MBGi_m(a,b,\bPhi,h)$, $i=1,2,3$ and $\bA$
is a constant $m$-dimensional nonsingular matrix. Then using
Theorem \ref{extended MBG}, the linear combination $\bA\bX_i\bA'$,
$i=1,2,3$ has the MBG distribution.
\end{remark}
\begin{remark}
In Theorem \ref{extended MBG}, it might be seemed that ``generalized noncentral" MBG distributions of
three types are defined.
\end{remark}

Entropy measures the uncertainty as confined in a distribution. Formally let
$(\chi ,\mathcal{B},\P )$ be a probability space, $f(.)$ is a density
function of matrix variate $\boldsymbol{X}$, associated with $\P $,
dominated by the $\sigma $-measure $\mu $ on $\chi $. The Shannon entropy
measures the expected information contained in the data and is equivalent to
the unpredicted component of a distribution. Then the well-known Shannon
entropy of $f$ is defined by
\begin{equation*}
E_{S}(f)=-\int_{\chi }f(\boldsymbol{X})\log f(\boldsymbol{X})\mathnormal{d}%
\mu .
\end{equation*}%
As an extension to the above measure, R\'{e}nyi entropy is defined as
\begin{equation*}
E_{R}(f)=\frac{1}{1-\nu }\log \int_{\chi }f^{\nu }(\boldsymbol{X})%
\mathnormal{d}\mu ,\quad \nu >0\ \mbox{and}\ \nu \neq 1.
\end{equation*}%
The additional parameter $\nu $, is used to describe complex behavior in
probability models and the associated process under study. R\'{e}nyi entropy
monotonically decreasing in $\nu $, while Shannon entropy is obtained from R%
\'{e}nyi for $\nu \uparrow 1$. For details see Zagrafos and
Nadarajah (2005). The R\'{e}nyi entropy for these distributions is
derived as follows.

\begin{theorem}~
\begin{enumerate}
\item[(i)] Let $\bX\sim MBG1_m(a,b,h,\bPhi)$. Then the R\'enyi entropy is
given by
\begin{eqnarray*}
E_R(f)&=&\frac{1}{1-\nu}\log\bigg[\left(\zeta_{a,b}^{(1)}\right)^\nu\sum_{t=0}^\infty\frac{u^{(t)}(0)}{t!}\sum_\tau\cr
&& \frac{\Gamma_m\left(\nu
a-\frac{1}{2}(\nu-1)(m+1),\tau\right)\Gamma_m\left(b
\nu-\frac{1}{2}(\nu-1)(m+1)\right)}{\Gamma_m\left(\nu a+\nu
b-(\nu-1)(m+1),\tau\right)}C_\tau(\bPhi)\bigg],
\end{eqnarray*}
\item[(ii)] Let $\bX\sim MBG2_m(a,b,h,\bPhi)$. Then the R\'enyi entropy is
given by
\begin{eqnarray*}
E_R(f)&=&\frac{1}{1-\nu}\log\bigg[\left(\zeta_{a,b}^{(2)}\right)^\nu\sum_{t=0}^\infty\frac{u^{(t)}(0)}{t!}\sum_\tau\cr
&& \frac{\Gamma_m\left(\nu
a-\frac{1}{2}(\nu-1)(m+1),\tau\right)\Gamma_m\left(b
\nu-\frac{1}{2}(\nu-1)(m+1),-\tau\right)}{\Gamma_m\left(\nu a+\nu
b-(\nu-1)(m+1)\right)}C_\tau(\bPhi)\bigg],
\end{eqnarray*}
\item[(iii)] Let $\bX\sim MBG3_m(a,b,h,\bPhi)$. Then the R\'enyi entropy is
given by
\begin{eqnarray*}
E_R(f)&=&\frac{1}{1-\nu}\log\bigg[\left(\zeta_{a,b}^{(3)}\right)^\nu\sum_{\tau,\kappa,\phi}\frac{(-1)^k(\nu
a+\nu b-(\nu-1)(m+1))_\kappa h^{(t)}(0)}{t!k!}\cr &&
\frac{\Gamma_m\left(\nu
a-\frac{1}{2}(\nu-1)(m+1),\phi\right)\Gamma_m\left(b
\nu-\frac{1}{2}(\nu-1)(m+1)\right)}{\Gamma_m(\nu a+\nu
b-(\nu-1)(m+1),\phi)}\frac{\left(\theta_\phi^{\kappa,\tau}\right)^2C_\phi(\bI_m)}{C_\tau(\bI_m)}C_\tau(\bPhi)\bigg],
\end{eqnarray*}
\end{enumerate}
where $u^{(t)}$ is the $t$-th derivative of $h^\nu$.
\end{theorem}\noindent \textbf{Proof:} By Definition \ref{MBG}, for MBG1
distribution, we have that
\begin{eqnarray*}
I(\nu )&=&\int_{\mathcal{I}_{m}}f^{\nu }(\boldsymbol{X})\mathnormal{d}%
\boldsymbol{X}\cr
&=&\left( \zeta _{a,b}^{(1)}\right) ^{\nu }\int_{\mathcal{I}%
_{m}}\det (\boldsymbol{X})^{\nu a-\frac{1}{2}(\nu -1)(m+1)-\frac{1}{2}%
(m+1)}\det (\boldsymbol{I}_{m}-\boldsymbol{X})^{\nu b-\frac{1}{2}(\nu
-1)(m+1)-\frac{1}{2}(m+1)}\cr\times h^{\nu }(\mathop{\rm tr}\nolimits(%
\boldsymbol{\Phi }\boldsymbol{X}))\mathnormal{d}\boldsymbol{X}.
\end{eqnarray*}%
Since $u(\mathop{\rm tr}\nolimits(\boldsymbol{\Phi }\boldsymbol{X}))=h^{\nu
}(\mathop{\rm tr}\nolimits(\boldsymbol{\Phi }\boldsymbol{X}))$ is a Borel
measurable function that admits a Taylor's series expansion in zonal
polynomials under some mild conditions, we get
\begin{eqnarray*}
I(\nu )&=&\left( \zeta _{a,b}^{(1)}\right) ^{\nu }\int_{\mathcal{I}_{m}}\det (%
\boldsymbol{X})^{\nu a-\frac{1}{2}(\nu -1)(m+1)-\frac{1}{2}(m+1)}\det (%
\boldsymbol{I}_{m}-\boldsymbol{X})^{\nu b-\frac{1}{2}(\nu -1)(m+1)-\frac{1}{2%
}(m+1)}\cr
&&\times u(\mathop{\rm tr}\nolimits(\boldsymbol{\Phi }\boldsymbol{X}%
))\mathnormal{d}\boldsymbol{X}\cr &=&\left( \zeta
_{a,b}^{(1)}\right) ^{\nu }\sum_{t=0}^{\infty
}\frac{u^{(t)}(0)}{t!}\sum_{\tau }\frac{\Gamma
_{m}\left( \nu a-\frac{1}{2}(\nu -1)(m+1),\tau \right) \Gamma _{m}\left(\nu b-\frac{1}{2}(\nu -1)(m+1)\right)}{%
\Gamma _{m}\left( \nu a+\nu b-(\nu -1)(m+1),\tau \right) }C_{\tau }(%
\boldsymbol{\Phi }).
\end{eqnarray*}%
Taking the logarithm from the above result gives the R\'{e}nyi
entropy. The proof for the other two types is the same.\hfill
$\blacksquare $

As the final important property, the joint distribution of eigenvalues for
three type of MBG distributions will be given in the next theorem.
\begin{theorem}
Let $g_i(\bLambda)$ denote the joint density function of
eigenvalues $(\bLambda=\diag(\lambda_1,\ldots,\lambda_m))$ of
$\bX_i\sim MBGi_m(a,b,\bPhi,h)$. Then we have
\begin{eqnarray*}
g_1(\bLambda)&=&\frac{\pi^{\frac{1}{2}m^2}\zeta_{a,b}^{(1)}}{\Gamma_m\left(\frac{1}{2}m\right)}
\sum_{k=0}^\infty\sum_\kappa
\frac{h^{(k)}(0)}{k!}\frac{C_\kappa(\bPhi)}{C_\kappa(\bI_m)}\cr
&&\times \prod_{i<j}^m (\lambda_i-\lambda_j)\prod_{i=1}^m\left(
\lambda_i^{a-\frac{1}{2}(m+1)}(1-\lambda_i)^{b-\frac{1}{2}(m+1)}\right)C_\kappa(\bLambda),\;1>\lambda_1>\ldots>\lambda_m>0,\cr
g_2(\bLambda)&=&\frac{\pi^{\frac{1}{2}m^2}\zeta_{a,b}^{(2)}}{\Gamma_m\left(\frac{1}{2}m\right)}
\sum_{k=0}^\infty\sum_\kappa
\frac{h^{(k)}(0)}{k!}\frac{C_\kappa(\bPhi)}{C_\kappa(\bI_m)}\cr
&&\times \prod_{i<j}^m (\lambda_i-\lambda_j)\prod_{i=1}^m\left(
\lambda_i^{a-\frac{1}{2}(m+1)}(1+\lambda_i)^{-(a+b)}\right)C_\kappa(\bLambda),\;\lambda_1>\ldots>\lambda_m>0,\cr
g_3(\bLambda)&=&\frac{\pi^{\frac{1}{2}m^2}\zeta_{a,b}^{(3)}}{\Gamma_m\left(\frac{1}{2}m\right)}
\sum_{k=0}^\infty\sum_\kappa
\frac{h^{(k)}(0)}{k!}\frac{C_\kappa(\bPhi)}{C_\kappa(\bI_m)}\cr
&&\times \prod_{i<j}^m (\lambda_i-\lambda_j)\prod_{i=1}^m\left(
\lambda_i^{a-\frac{1}{2}(m+1)}(1-\lambda_i)^{b-\frac{1}{2}(m+1)}(1+\lambda_i)^{-(a+b)}\right)C_\kappa(\bLambda),\;1>\lambda_1>\ldots>\lambda_m>0.\cr
\end{eqnarray*}
\end{theorem}
\noindent\textbf{Proof:} From Theorem 3.2.17 of Muirhead (2005), the density
of $\boldsymbol{\Lambda}$, for MBG1 distribution, is given by
\begin{eqnarray*}
g_1(\boldsymbol{\Lambda}) &=&\frac{\pi^{\frac{1}{2}m^2}}{\Gamma_m\left(\frac{%
1}{2}m\right)}\prod_{i<j}^m (\lambda_i-\lambda_j)\int_{\O (m)}f(\boldsymbol{H%
}\boldsymbol{\Lambda}\boldsymbol{H}^{\prime })\mathnormal{d}\boldsymbol{H}%
\cr &=&\frac{\pi^{\frac{1}{2}m^2}}{\Gamma_m\left(\frac{1}{2}m\right)}%
\prod_{i<j}^m (\lambda_i-\lambda_j)\zeta_{a,b}^{(1)}\det(\boldsymbol{\Lambda}%
)^{a-\frac{1}{2}(m+1)}\det(\boldsymbol{I}_m-\boldsymbol{\Lambda})^{b-\frac{1%
}{2}(m+1)}\cr &&\times\int_{\O (m)}h\left(\mathop{\rm tr}\nolimits%
\boldsymbol{\Phi}\boldsymbol{H}\boldsymbol{\Lambda}\boldsymbol{H}^{\prime
}\right)\mathnormal{d}\boldsymbol{H}.
\end{eqnarray*}
Making use of Eq. (36) of Muirhead (2005) follows
\begin{eqnarray*}
\int_{\O (m)}h\left(\mathop{\rm tr}\nolimits\boldsymbol{\Phi}\boldsymbol{H}%
\boldsymbol{\Lambda}\boldsymbol{H}^{\prime }\right)\mathnormal{d}\boldsymbol{%
H}&=&\sum_{k=0}^\infty\sum_\kappa \frac{h^{(k)}(0)}{k!}\int_{\O %
(m)}C_\kappa\left(\boldsymbol{\Phi}\boldsymbol{H}\boldsymbol{\Lambda}%
\boldsymbol{H}^{\prime }\right)\mathnormal{d}\boldsymbol{H}\cr %
&=&\sum_{k=0}^\infty\sum_\kappa \frac{h^{(k)}(0)}{k!}\frac{C_\kappa(%
\boldsymbol{\Phi})C_\kappa(\boldsymbol{\Lambda})}{C_\kappa(\boldsymbol{I}_m)}%
.
\end{eqnarray*}
Thus the final result immediately follows by noting that $\det(\boldsymbol{%
\Lambda})^{a-\frac{1}{2}(m+1)}=\prod_{i=1}^m \lambda_i^{a-\frac{1}{2}(m+1)}$
and $\det(\boldsymbol{I}_m-\boldsymbol{\Lambda})^{b-\frac{1}{2}%
(m+1)}=\prod_{i=1}^m (1-\lambda_i)^{b-\frac{1}{2}(m+1)}$. The
proof for the other two types can be achieved in a similar
fashion.\hfill$\blacksquare$

\section{Discussion}

Some research questions that emanates from Definition \ref{definnition
kernel oriented}\ are highlighted below.

\textbf{4.1)} Different well-known matrix variate distributions,
as well as new case(s) follows from this Definition \ref{MBG}. In
2002, Nagar and Gupta proposed the matrix variate Kummer beta
distribution extending the work of Ng and Kotz (1995). Now in this
section we will focus on the matrix variate Kummer beta (MKB1/2/3)
distribution as special case of MBG1/2/3 distribution, that can be
obtained by taking $h(x)=\exp (-x)$ in Definition \ref{MBG}. The
first two types are well-known in literature, however MKB type 3
is new. In this regard, we have the following general definition.

\begin{definition}\label{MKB}
Let $\Re(a)>(m-1)/2$, $\Re(b)>(m-1)/2$ and $\bPhi\in\S_m$. Then
the random symmetric matrix $\bX$ of dimension $m$ is said to have
\begin{enumerate}
\item[(i)] MKB1 distribution with parameters
$a$, $b$ and $\bPhi$ and shape generator $h$, denoted by $\bX\sim
MKB1_m(a,b,\bPhi)$, if it has the following density function
\begin{eqnarray*}
f(\bX)=K_{a,b}^{(1)}\det(\bX)^{a-\frac{1}{2}(m+1)}\det(\bI_m-\bX)^{b-\frac{1}{2}(m+1)}\exp(-\tr(\bPhi\bX)),
\quad \bX\in\I_m,
\end{eqnarray*}
where $\left(K_{a,b}^{(1)}\right)^{-1}=B_m(a,b) \;
_1F_1(a,a+b;-\bPhi)$ (see Nagar and Gupta, 2002).
\item[(ii)] MKB2 distribution with parameters
$a$, $b$ and $\bPhi$ and shape generator $h$, denoted by $\bX\sim
MKB2_m(a,b,\bPhi)$, if it has the following density function
\begin{eqnarray*}
f(\bX)=K_{a,b}^{(2)}\det(\bX)^{a-\frac{1}{2}(m+1)}\det(\bI_m+\bX)^{-(a+b)}\exp(-\tr(\bPhi\bX)),
\quad \bX\in\S_m,
\end{eqnarray*}
where using Lemma 5 of Khatri (1966)
\begin{eqnarray*}
\left(K_{a,b}^{(2)}\right)^{-1}&=&\sum_{k=0}^\infty\sum_\kappa\frac{1}{k!}\int_{\S_m}\det(\bX)^{a-\frac{1}{2}(m+1)}\det(\bI_m+\bX)^{-(a+b)}
C_\kappa(-\bPhi\bX)\textnormal{d}\bX\cr
&=&\sum_{k=0}^\infty\sum_\kappa\frac{1}{k!}\frac{\Gamma_m(a,\kappa)\Gamma_m(b,-\kappa)}{\Gamma_m(a+b)}C_\kappa(-\bPhi).
\end{eqnarray*}
\item[(iii)] MKB3 distribution with parameters
$a$, $b$ and $\bPhi$ and shape generator $h$, denoted by $\bX\sim
MKB3_m(a,b,\bPhi)$, if it has the following density function
\begin{eqnarray*}
f(\bX)=K_{a,b}^{(3)}\det(\bX)^{a-\frac{1}{2}(m+1)}\det(\bI_m-\bX)^{b-\frac{1}{2}(m+1)}\det(\bI_m+\bX)^{-(a+b)}\exp(-\tr(\bPhi\bX)),
\quad \bX\in\I_m,
\end{eqnarray*}
where using Eq. (2.8) of Davis (1979) and Eq. (3.28) of Chikuse
(1980),
\begin{eqnarray*}
\left(K_{a,b}^{(3)}\right)^{-1}&=&\sum_{k=0}^\infty\sum_\kappa\frac{(-1)^k}{k!}\sum_{t=0}^\infty\sum_\tau\frac{(-1)^t(a+b)_\tau}{t!}\cr
&&
\int_{\S_m}\det(\bX)^{a-\frac{1}{2}(m+1)}\det(\bI_m-\bX)^{b-\frac{1}{2}(m+1)}
C_\kappa(\bPhi\bX)C_\tau(\bX)\textnormal{d}\bX\cr &=&
\sum_{\tau,\kappa,\phi}\frac{(-1)^{t+k}(a+b)_\tau
}{t!k!}
                          \frac{\Gamma_m(b)\Gamma_m(a,\phi)}{\Gamma_m(a+b,\phi)}\frac{\left(\theta_\phi^{\kappa,\tau}\right)^2
                          C_\phi(\bI_m)}{C_\kappa(\bI_m)}C_\kappa(\bPhi).
\end{eqnarray*}
\end{enumerate}
\end{definition}
To illustrate the effect of the shape structure ascribed to the
Borel measurable function combined with the kernel of a
statistical distribution, graphical representations are provided
for some cases. Some 3-dimensional graphical representations are
provided in Figure 1 for different parameter values of $\phi$;
2-dimensional representations are also given for different set
parameters $\boldsymbol{\Theta }=(a,b,\phi )$ in Figure 2, for
$h(x)=\exp(-x)$.

\begin{figure}[tbp]
\begin{center}
\begin{tabular}{ccc}
\includegraphics[scale=0.3]{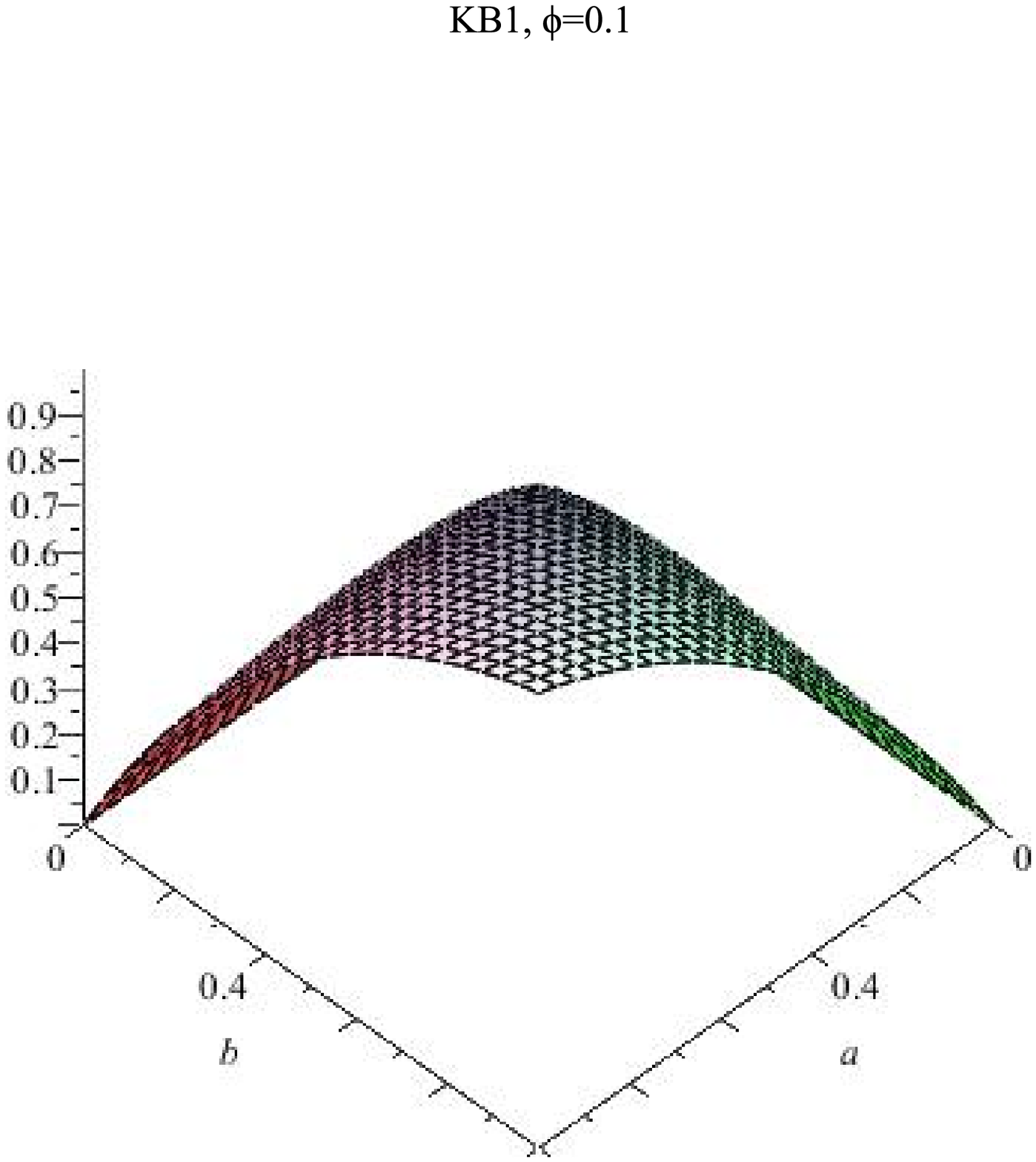} & %
\includegraphics[scale=0.3]{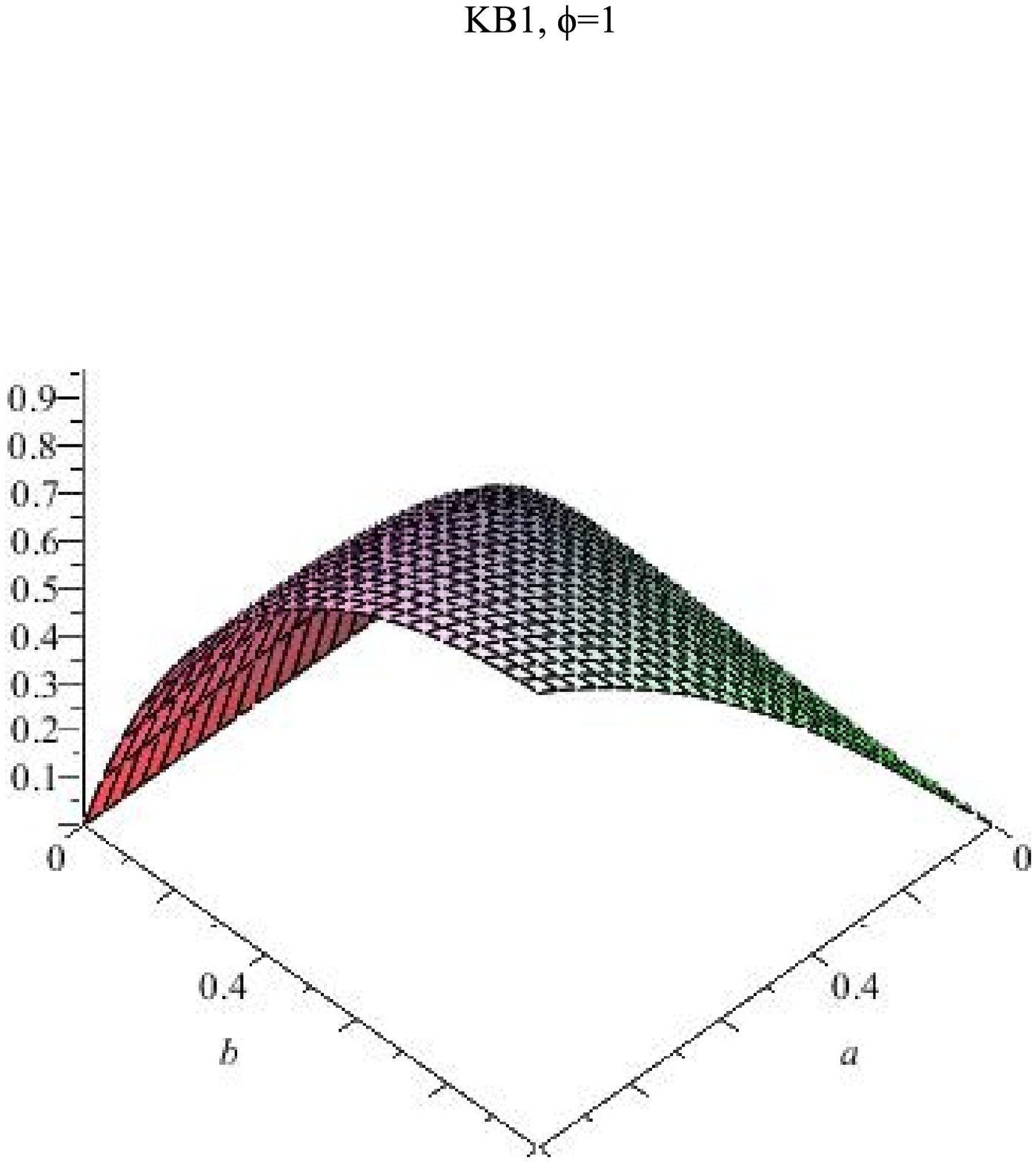} & %
\includegraphics[scale=0.3]{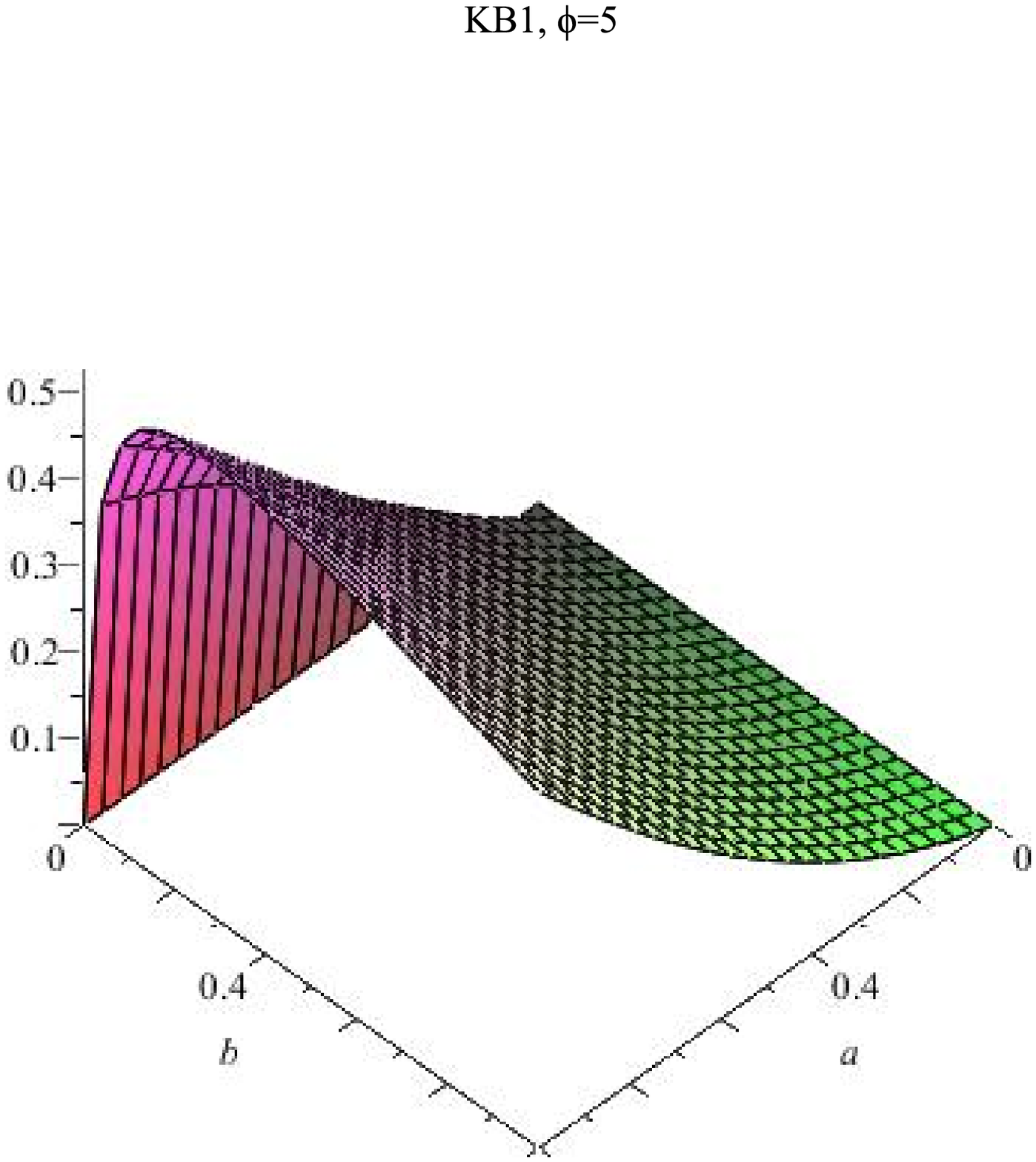} \\
\includegraphics[scale=0.3]{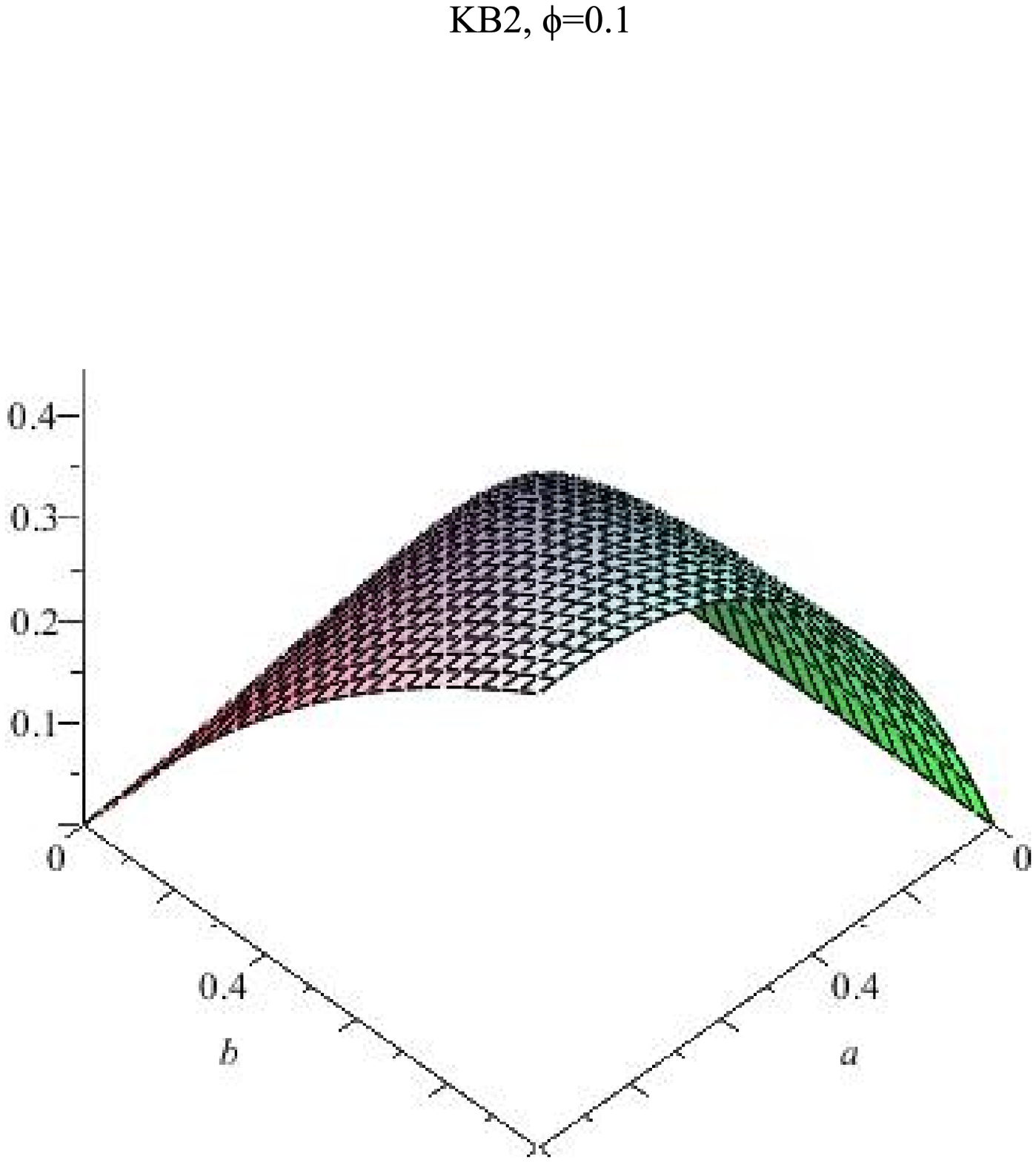} & %
\includegraphics[scale=0.3]{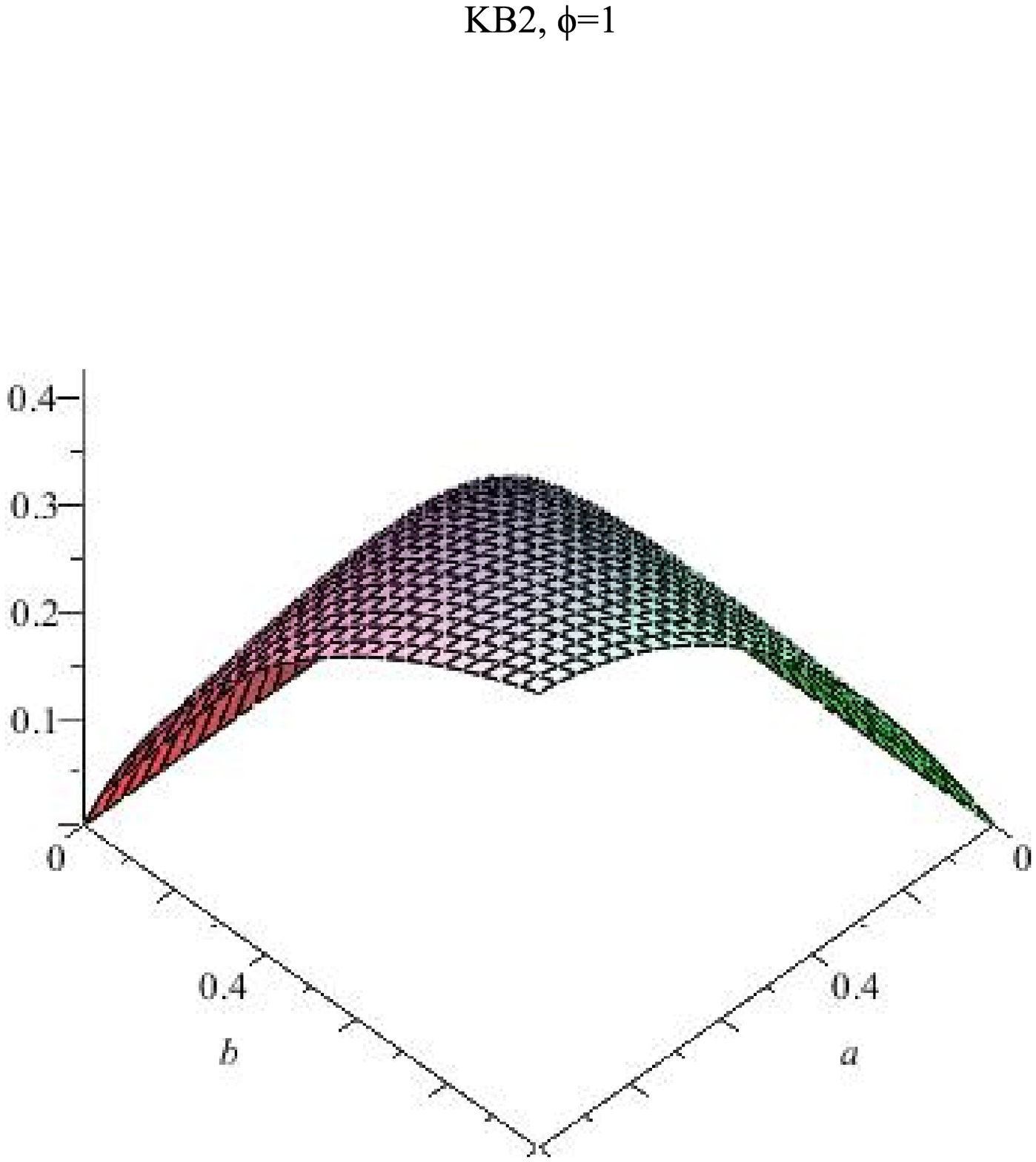} & %
\includegraphics[scale=0.3]{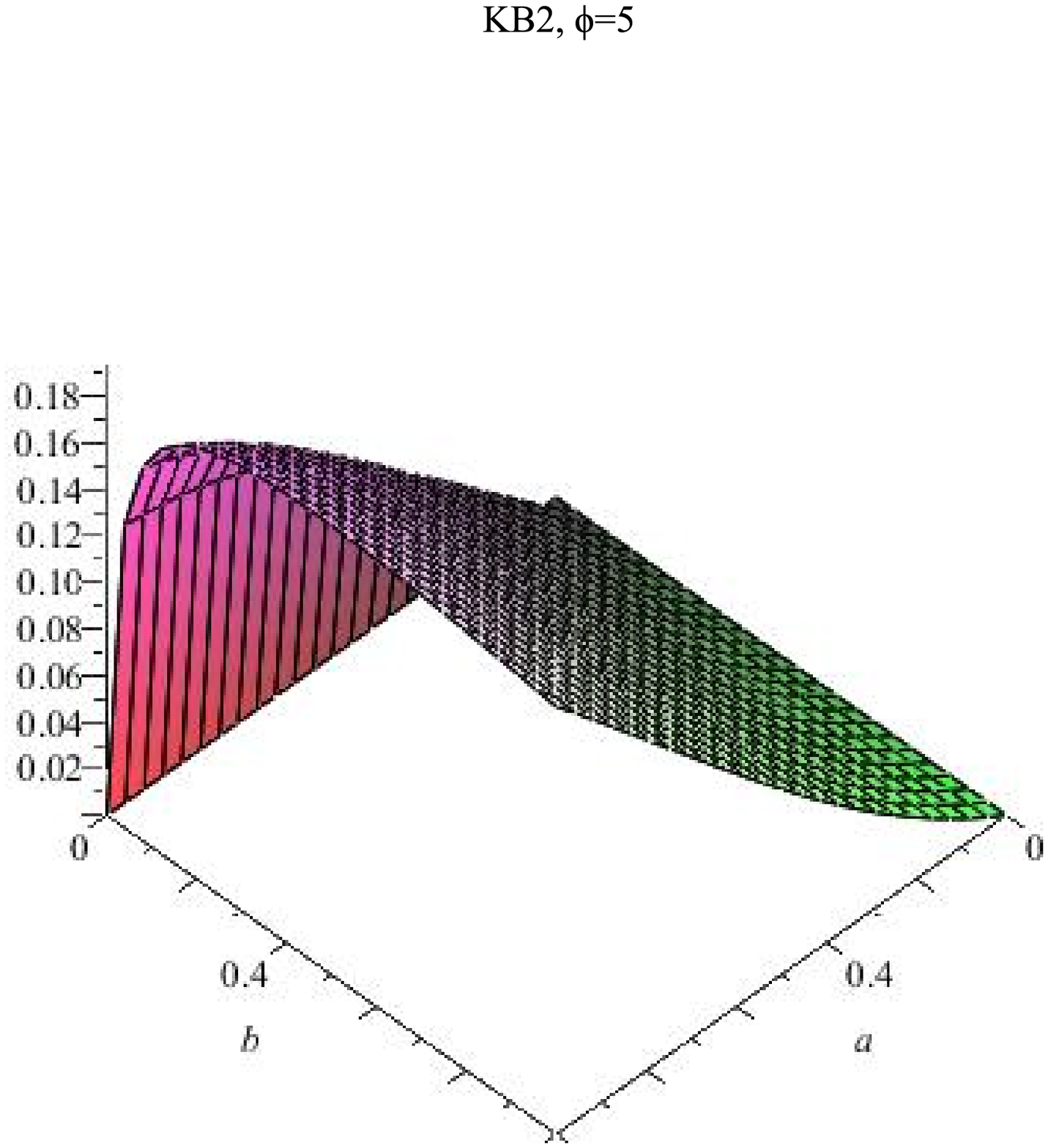} \\
\includegraphics[scale=0.3]{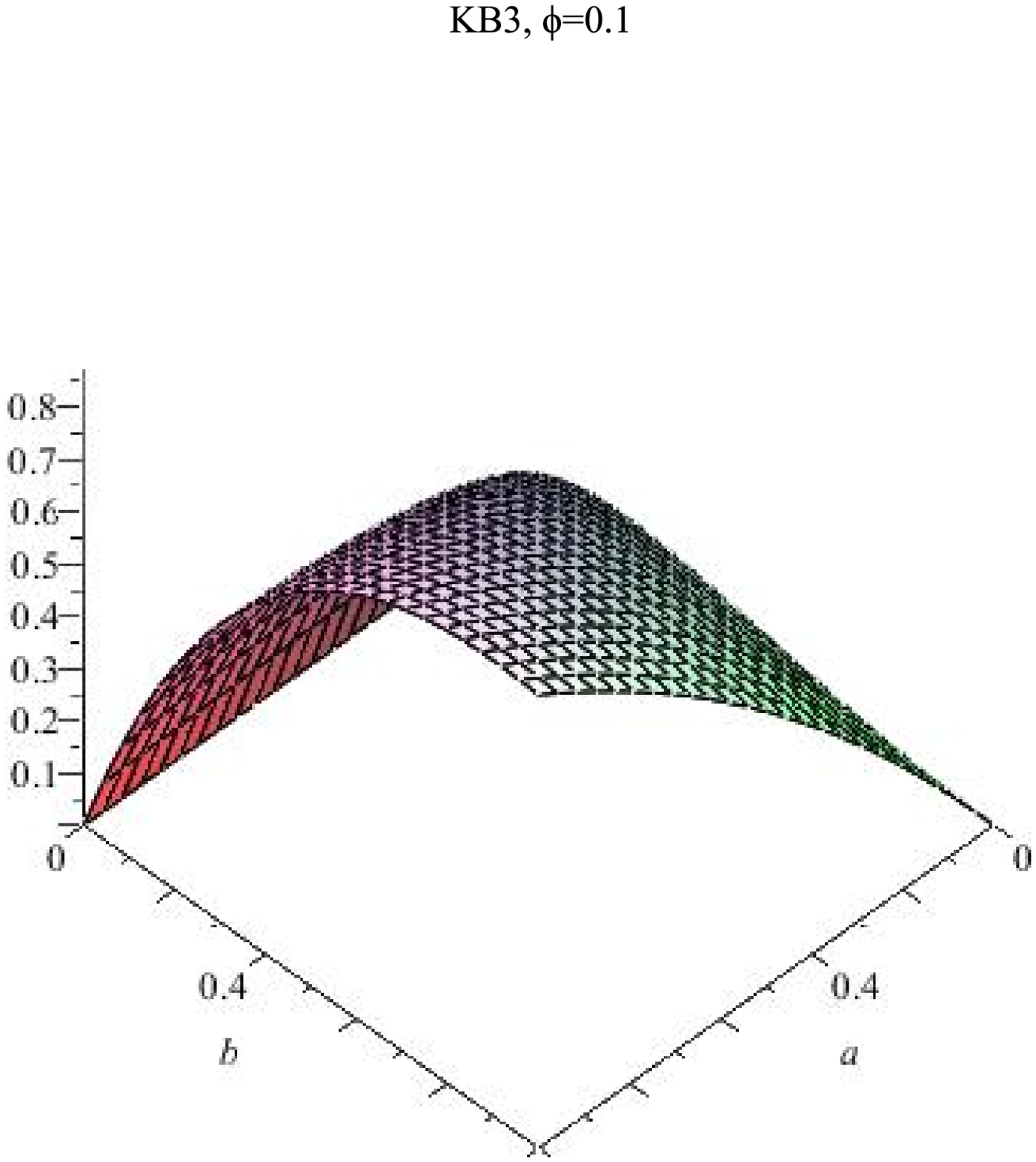} & %
\includegraphics[scale=0.3]{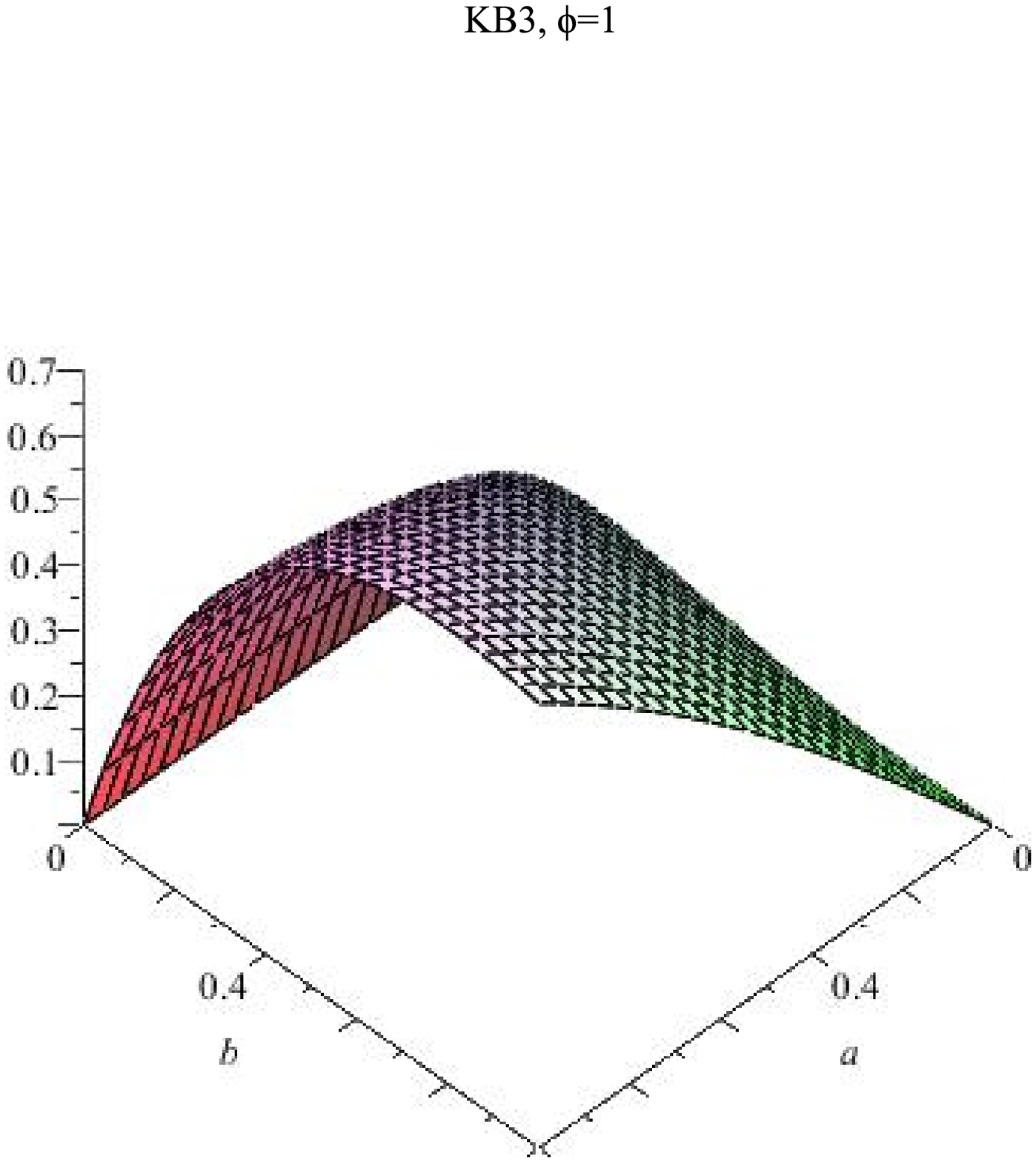} & %
\includegraphics[scale=0.3]{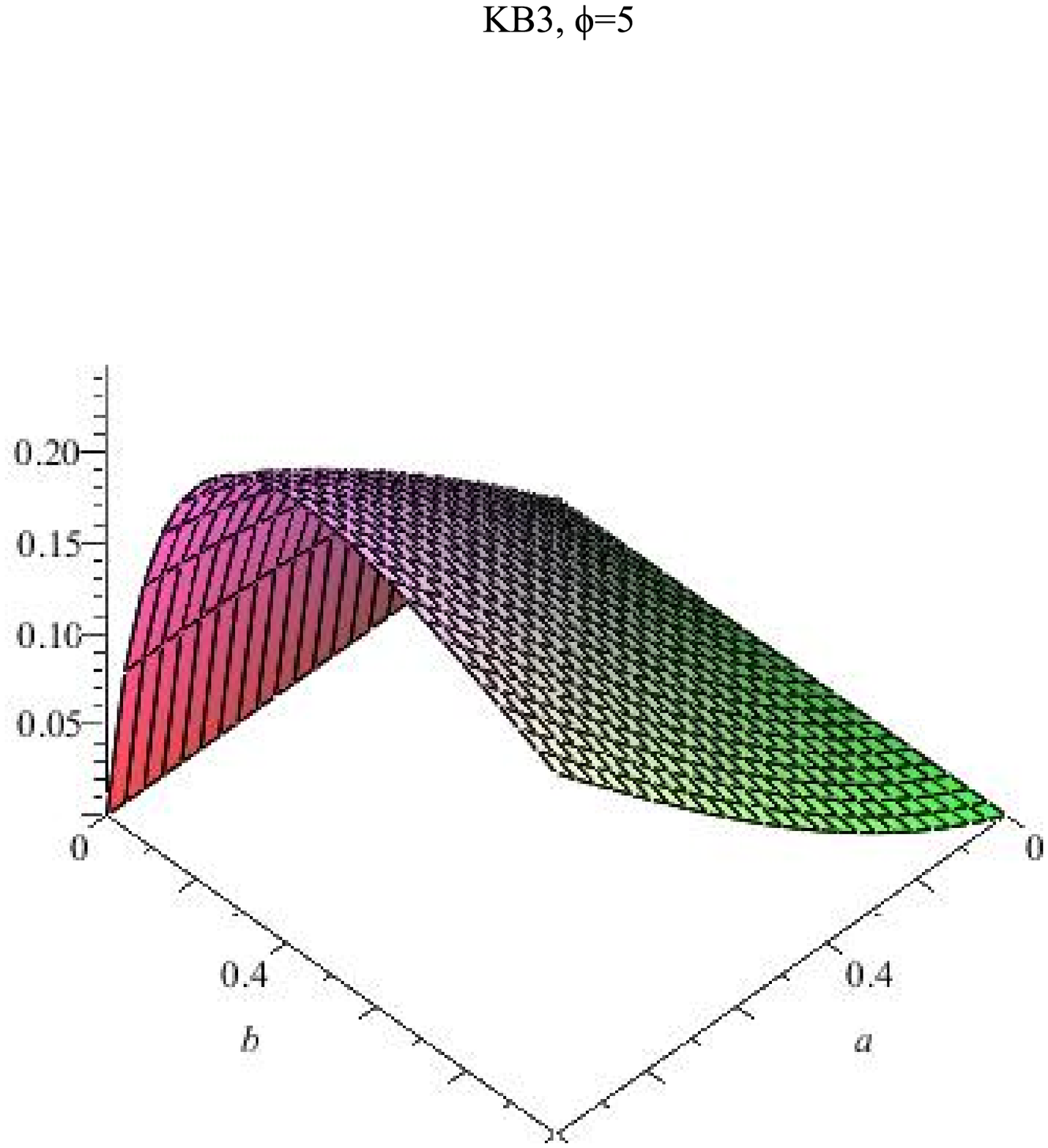} \\
&  &
\end{tabular}%
\end{center}
\caption{Kummer beta distributions for different $\protect\phi $
parameter values}
\end{figure}

\begin{figure}[tbp]
\begin{center}
\begin{tabular}{c}
\includegraphics[scale=0.4]{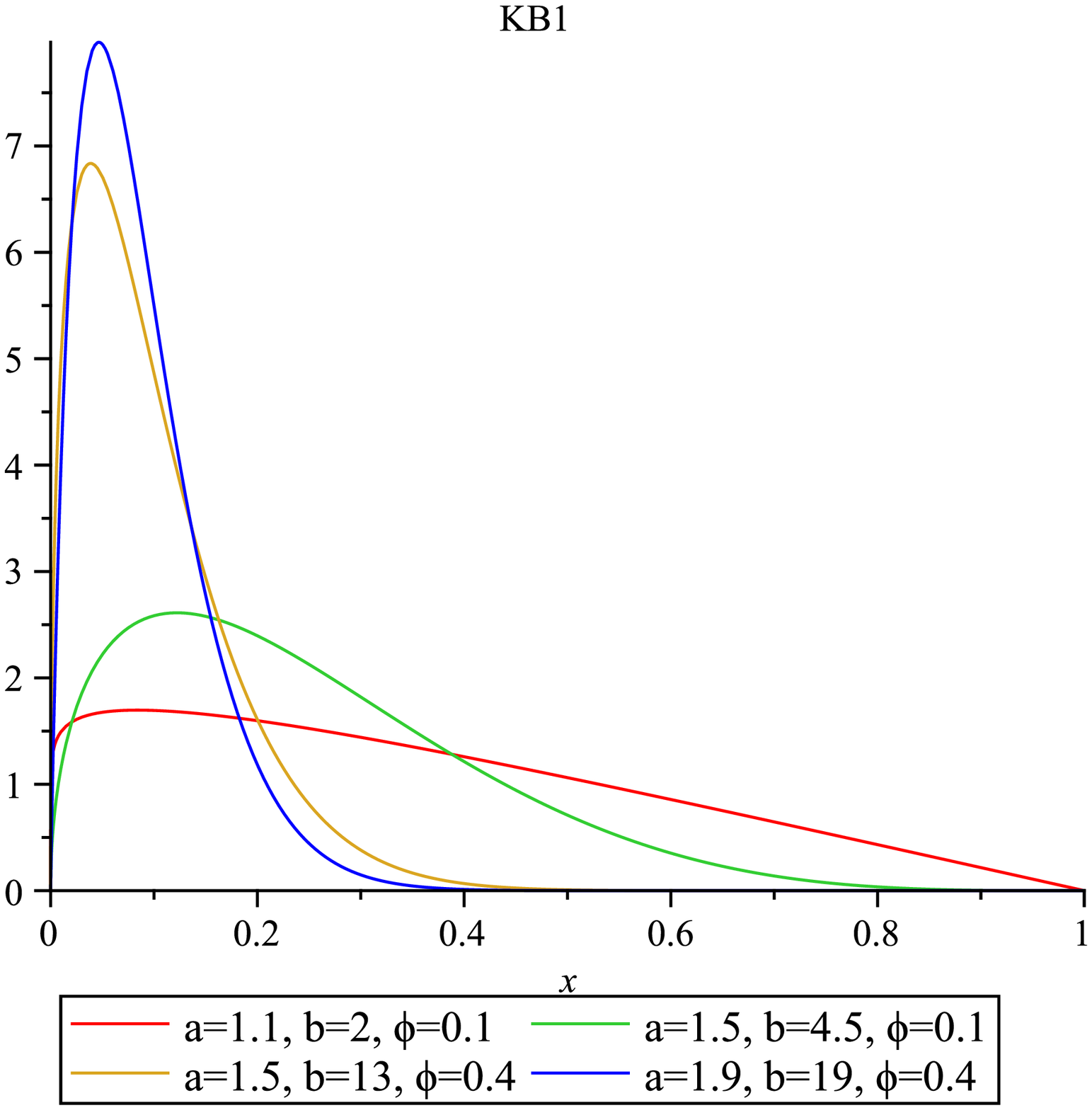} \\
\includegraphics[scale=0.4]{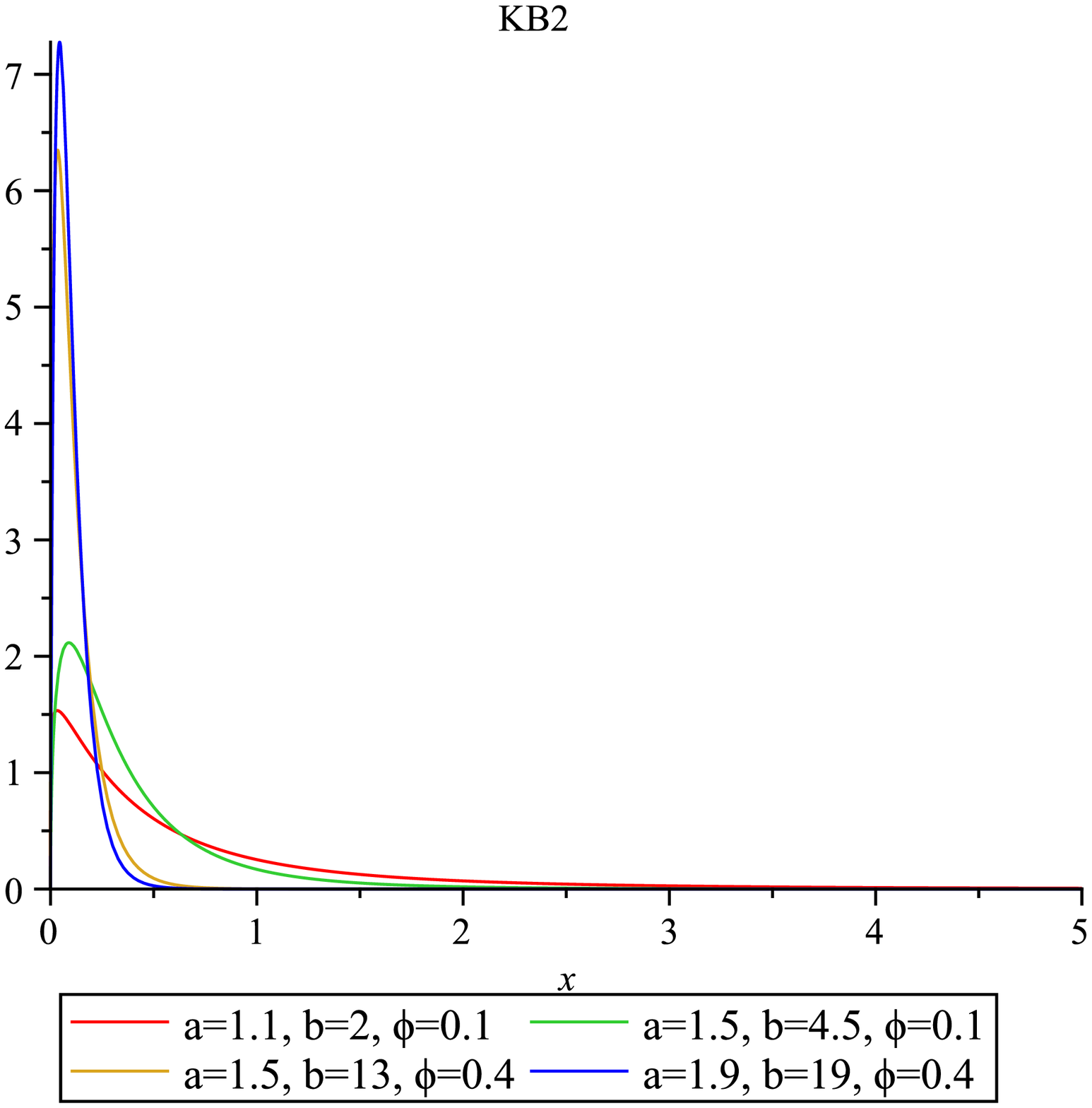} \\
\includegraphics[scale=0.4]{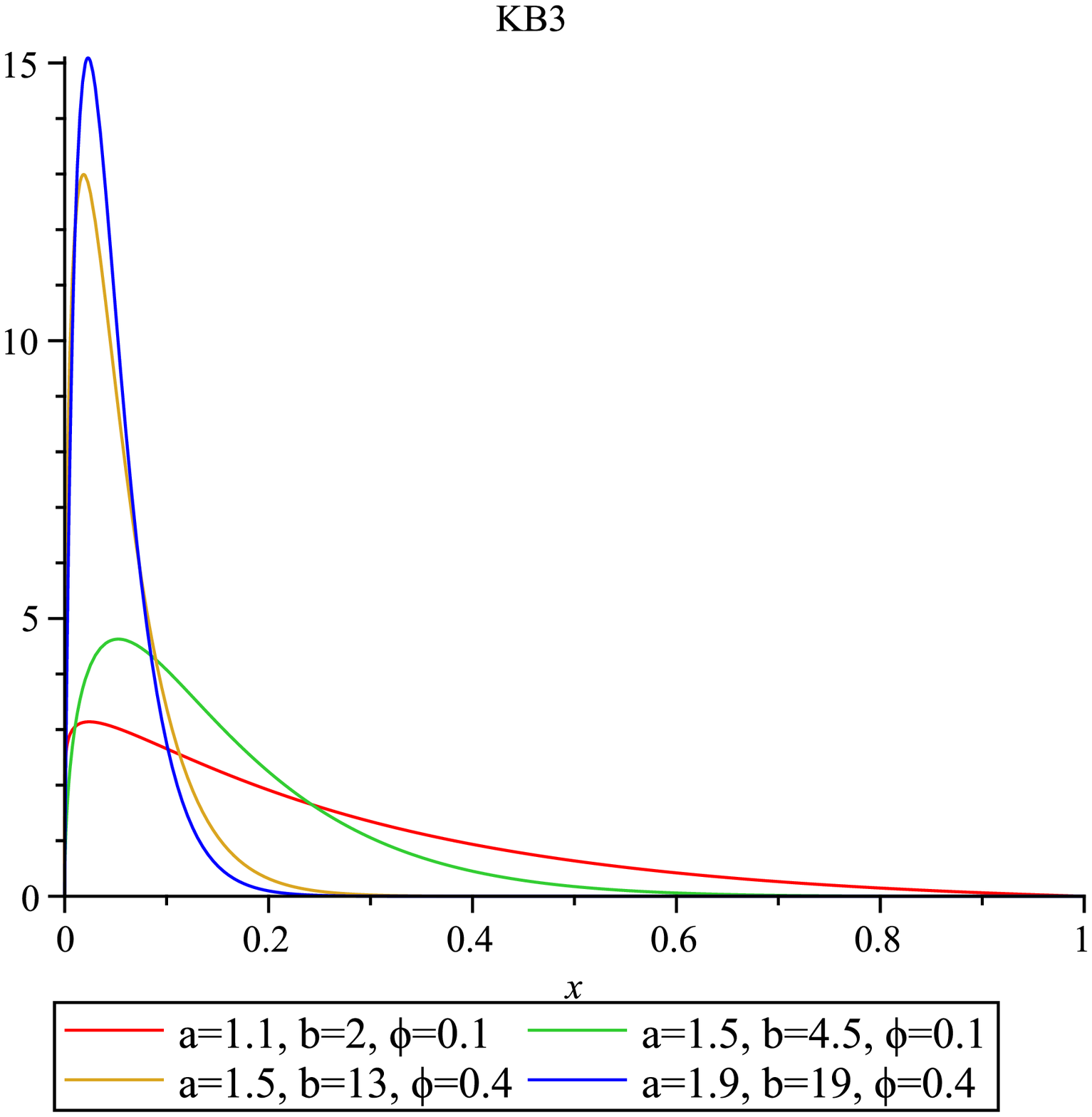} \\
\end{tabular}%
\end{center}
\caption{2-dimensional representation for univariate Kummer beta
distributions}
\end{figure}

\textbf{4.2)} According to Definition 1.1 all PK functions were
functions of the trace argument, however it is also possible to
extend the definition to include for example the PK function with
the determinant as argument. In this regard, we propose the
following definition.
\begin{definition}
The random symmetric matrix $\bX$ of dimension $m$ is said to have
matrix variate beta generator distribution with parameters $a$,
$b$ and $\bPhi$ and shape generator $h$, if it has the following
density
\begin{enumerate}
\item{(i)} 1st kind
\begin{eqnarray*}
f(\bX) &=& \left(\sum_{t=0}^\infty
\frac{h^{(t)}(0)B_m(a+t,b)}{t!}\det(\bPhi)^t\right)^{-1}\cr
&&\det(\bX)^{a-\frac{1}{2}(m+1)}\det(\bI_m-\bX)^{b-\frac{1}{2}(m+1)}\;
h(\det(\bPhi\bX)),\quad\bX\in\I_m,
\end{eqnarray*}
\item{(ii)} 2nd kind
\begin{eqnarray*}
f(\bX) &=& \left(\sum_{t=0}^\infty
\frac{h^{(t)}(0)B_m(a+t,b-t)}{t!}\det(\bPhi)^t\right)^{-1}\cr
&&\det(\bX)^{a-\frac{1}{2}(m+1)}\det(\bI_m+\bX)^{-(a+b)}\;
h(\det(\bPhi\bX)),\quad\bX\in\S_m,
\end{eqnarray*}
\item{(iii)} 3rd kind
\begin{eqnarray*}
f(\bX) &=& \left(\sum_{t=0}^\infty
\frac{h^{(t)}(0)B_m(a+t,b)}{2^{m(a+t)}t!}\det(\bPhi)^t\right)^{-1}\cr
&&\det(\bX)^{a-\frac{1}{2}(m+1)}\det(\bI_m-\bX)^{b-\frac{1}{2}(m+1)}\det(\bI_m+\bX)^{-(a+b)}\;
h(\det(\bPhi\bX)),\quad\bX\in\I_m,
\end{eqnarray*}
where $\Re(a),\Re(b)>(m-1)/2$, $\bPhi\in\I_m$ is a symmetric
complex matrix, $h(.)$ is a Borel measurable function that admits
a Taylor series expansion in zonal polynomials
\end{enumerate}
\end{definition}

\textbf{4.3)} It is known that the Wilks' statistic plays the same
role in multivariate analysis as the F statistic plays in
univariate analysis. Bekker et al. (2011) derived an exact
expression for the non-null distribution of the Wilks' statistics.
Bekker et al. (2012) proposed new multivariate test statistics and
their exact distributions. Thus it is of interest to find the
distribution of the determinant where the matrix variate has the
MGB$i(i=1,2,3)$ distribution leading to generalized Wilks'
statistics.

\begin{theorem}
Let $\bX_i\sim MBGi_m(a,b,\bPhi,h)$, $i=1,2,3$. Then
\begin{enumerate}
\item[1.] $y=\det(\bX_1)$ has the following density function
\begin{equation*}
f(y)=\zeta_{a,b}^{(1)}\Gamma_m(b)\sum_{t=0}^\infty\frac{h^{(t)}(0)}{t!}\sum_\tau
C_\tau(\bPhi) G_{m,m}^{m,0}\left(y\big|\begin{array}{c}
c_1,\ldots,c_m\\d_1,\ldots,d_m\end{array}\right),
\end{equation*}
where $c_j=a+b+t_j-\frac{1}{2}(j+1)$ and
$d_j=a+t_j-\frac{1}{2}(j+1)$.
\item[2.] $y=\det(\bX_2)$ has the following density function
\begin{equation*}
f(y)=\zeta_{a,b}^{(2)}\sum_{t=0}^\infty\frac{h^{(t)}(0)}{t!}\sum_\tau
\Gamma_m(b,-\tau)C_\tau(\bPhi)
G_{m,m}^{m,0}\left(y\big|\begin{array}{c}
c_1,\ldots,c_m\\d_1,\ldots,d_m\end{array}\right),
\end{equation*}
where $c_j=a+b+t_j-\frac{1}{2}(j+1)$ and
$d_j=a+t_j-\frac{1}{2}(j+1)$.
\item[3.] $y=\det(\bX_3)$ has the following density function
\begin{equation*}
f(y)=\zeta_{a,b}^{(3)}\Gamma_m(b)\sum_{\tau,\kappa,\phi}\frac{(-1)^kh^{(t)}(0)\left(\theta_\phi^{\kappa,\tau}\right)^2}{t!k!C_\tau(\bI_m)}
C_\tau(\bPhi) G_{2m,2m}^{2m,0}\left(y\big|\begin{array}{c}
c_1,\ldots,c_{2m}\\d_1,\ldots,d_{2m}\end{array}\right),
\end{equation*}
where
\begin{eqnarray*}
c_j&=&\left\{\begin{array}{cc}
a+b-1+\phi_{\frac{i+1}{2}}-\frac{1}{4}(i-1), &
i=1,3,5,\ldots,2m-1\\
a+b-1-\frac{1}{4}(i-2), & i=2,6,10,\ldots,2m\end{array}\right.\cr
d_j&=&\left\{\begin{array}{cc}
a+b-1+k_{\frac{i+1}{2}}-\frac{1}{4}(i-1), &
i=1,3,5,\ldots,2m-1\\
a-1+\phi_{\frac{i}{2}}-\frac{1}{4}(i-2), &
i=2,6,10,\ldots,2m\end{array}\right.
\end{eqnarray*}
\end{enumerate}
where $G(.)$ denotes the Meijer's G function.
\end{theorem}
\noindent\textbf{Proof:} We only give the proof for item 1; the
proofs of other two types are the same. Using Theorem
\ref{expectation}, the Mellin transform is given by (see Mathai,
1993)
\begin{eqnarray*}
\mathcal{M}_f(r)&\equiv& E\left(\det(\boldsymbol{X})^{r-1}\right)\cr %
&=&\zeta_{a,b}^{(1)}\sum_{t=0}^\infty \frac{h^{(t)}(0)}{t!}\sum_\tau\frac{%
\Gamma_m(a+r-1,\tau)\Gamma_m(b)}{\Gamma_m(a+b+r-1,\tau)}C_\tau(\boldsymbol{\Phi}%
).
\end{eqnarray*}
Thus the distribution of $y$ is uniquely obtained from the inverse
Mellin transform of the above and the definition of the Meijer's
G-function, $G(.)$. The proof is complete.\hfill$\blacksquare$

In this paper we developed the conventional matrix variate beta
distributions to more general ones, where the \textit{kernel of a
matrix variate beta type 1/2/3 with a further Borel measurable
function} of scalar value were combined. Important statistical
characteristics were derived such as the moment generating
function as well as the joint density function of eigenvalues. The
matrix variate Kummer beta distributions were discussed as special
cases. The authors are currently developing more theory and
results based on the principle of Definition 1.1 and the theory
applied in the paper. The program of work will include amongst
others the noncentral beta as kernel combined with numerous
different generators (see Arashi et al., 2013 and Van Niekerk et
al., 2013).

With the similar idea of generating new families of matrix-variate distributions,
another families of distributions can be generated by utilizing the ``Wishart-type kernels" combined
with an unknown Borel measurable function of trace and/or determinant operators. In this case the
algebra will need evaluating "Laplace-type integrals" involving zonal polynomials. (Refer to Bekker et al., 2013)

We deem that the proposed results in this paper should
stimulate research and applications beyond the known matrix
variate distributions.

\section*{Acknowledgments}

The authors would like to hereby acknowledge the support of the StatDisT
group. This work is based upon research supported by the National Research
foundation, South Africa (Incentive Funding for Rated Researchers) and VC
Post-Doctoral Fellowship of the University of Pretoria.

\section*{References}

\baselineskip=12pt

\noindent\hangindent 25pt M. Arashi, A. Bekker, J. J. J. Roux, and
J. Van Niekerk, (2013). Noncentral beta kernel oriented generator
family. \emph{Technical Report}, University of Pretoria, South
Africa.

\noindent\hangindent 25pt A. Bekker, J. J. J. Roux, and M. Arashi,
(2011). Exact nonnull distribution of Wilks' statistic: The ratio
and product of independent components, \emph{J. Mult. Anal.}, 102,
619-628.

\noindent\hangindent 25pt A. Bekker, J. J. J. Roux, R. Ehlers and
M. Arashi, (2012). Distribution of the product of determinants of
noncentral bimatrix beta variates, \emph{J. Mult. Anal.}, 109,
73-87.

\noindent\hangindent 25pt A. Bekker, M. Arashi, J. J. J. Roux, and J. van Niekerk, (2013). Wishart kernel
oriented generator distribution: Weighted Wishart distribution, Technical Report ISBN: 978-1-
77592-067-0, University of Pretoria, South Africa.

\noindent\hangindent 25pt B. W. Brown, M. S. Floyd, L. B. Levy,
(2002). The log F: a distribution for all seasons. \emph{Comput.
Statist.}, 17, 47-58.

\noindent \hangindent25pt Yasuko Chikuse, (1980). Invariant
polynomials with matrix arguments and their applications in
\emph{Multivariate Statistical Analysis}, 1, 54-68.

\noindent\hangindent 25pt A. W. Davis, (1979). Invariant
polynomials with two matrix arguments extending the zonal
polynomials: Applications to multivariate distribution theory,
\emph{Ann. Inst. Statist. Math.}, 31(A), 465-485.

\noindent\hangindent 25pt A. W. Davis, (1980). Invariant
polynomials with two matrix arguments, extending the zonal
polynomials, \emph{Multivariate Ananlysis-V} (ed. P. R.
Krishnaiah), 287-299.

\noindent\hangindent 25pt R. Ehlers, (2011). \emph{Bimatrix
Variate Distributions of Wishart Ratios With Application},
Unpublished PhD Dissertation, University of Pretoria.

\noindent\hangindent 25pt M. Ferreia, M. I. Gomez, and V. Leiva,
(2012). On an extreme value version of the Birnbaum-Sanders distribution, \emph{REVSTAT}%
, 10(2), 181-210.

\noindent\hangindent 25pt Arjun K. Gupta and Daya K. Nagar, (2000a). \emph{%
Matrix variate distributions}, Chapman and Hall / CRC, Boca Raton.

\noindent\hangindent 25pt Arjun K. Gupta and Daya K. Nagar,
(2000b). Matrix-variate beta distribution, \emph{Int. J. Math.
Sci.}, 24(7), 49-459.

\noindent\hangindent 25pt Arjun K. Gupta and Daya K. Nagar,
(2002). Matrix-variate Kummer-beta distribution, \emph{J. Aust.
Math. Sci.}, 73, 11-25.

\noindent\hangindent 25pt Arjun K. Gupta and Daya K. Nagar,
(2006). A
Generalized Matrix Variate Beta Distribution, \emph{Int. J. Appl. Math. Sci.}%
, 3(1), 21-36.

\noindent\hangindent 25pt Arjun K. Gupta and Daya K. Nagar,
(2009). Properties of matrix variate beta type 3 distribution,
\emph{Int. J. Mathematics Math. Sci.},
http://dx.doi.org/10.1155/2009/308518.

\noindent\hangindent 25pt A. T. James, (1961). Zonal polynomials
of the real positive definite symmetric matrices, \emph{Ann.
Math.}, 35, 456-469.

\noindent\hangindent 25pt A. T. James, (1964). Distribution of
matrix variate and latent roots derived from normal samples,
\emph{Ann. Math. Statist.}, 35, 475-501.

\noindent\hangindent 25pt M. C. Jones, (2004). Family of
distributions arising from distribution of order statistics.
\emph{Test} 13, 1-43.

\noindent\hangindent 25pt C. G. Khatri, (1966). On certain
distribution problemsbased on positive definite quadratic
functions in normal vector, \emph{Ann. Math. Statist.}, 37,
468-479.

\noindent\hangindent 25pt P. Kvam, (2008). Length bias in the
measurements of carbon nanotubes, \emph{Technometrics}, 50(4),
462-467.

\noindent\hangindent 25pt A. M. Mathai, (1993). \emph{A Handbook
of Generalized Special Functions for Statistical and Physical
Sciences}, Clarendon Press, Oxford.

\noindent\hangindent 25pt Rob J. Muirhead, (2005). \emph{Aspects
of Multivariate Statistical Theory}, 2nd Ed., John Wiley, New
York.

\noindent\hangindent 25pt S. Nadarajah and S. Kotz, (2004). The
beta Gumbel distribution, \emph{Math. Prob. Eng.}, 10, 323-332.

\noindent\hangindent 25pt S. Nadarajah and S. Kotz, (2006). The
beta exponential distribution, \emph{Reliab. Eng. Syst. Saf.}, 91,
689-697.

\noindent\hangindent 25pt S. Nadarajah and S. Kotz, (2006). Some
beta distributions, \emph{Bull. Brazilian Math. Soc.}, New Series,
37(1), 103-125.

\noindent\hangindent 25pt Daya K. Nagar and Arjun K. Gupta,
(2002). Matrix-variate Kummer-beta distribution, \emph{J.
Australian Math. Soc.}, 73(1), 11-25.

\noindent\hangindent 25pt Daya K. Nagar, Alejandro Rold\'a-Correa
and Arjun K. Gupta, (2013). Extended matrix variate gamma and beta
functions, \emph{J. Mult. Anal.}, 122, 53-69.

\noindent\hangindent 25pt A. K. Nanda and K. Jain (1999). Some
weighted distribution results on univariate and bivariate cases,
\emph{J. Statist. Plann. Inf.}, 77, 169-180.

\noindent\hangindent 25pt J. Navarro, J. M. Ruiz and Y. DeL
Aguila, (2006). Multivariate weighted distributions: a review and
some extensions, \emph{Statistics}, 40(1), 51-64.

\noindent\hangindent 25pt K. W. Ng and S. Kotz, (1995). Kummer
gamma and Kummer beta univariate and multivariate distributions,
Research Report No. 84, (Department of Statistics, The University
of Hong Kong, Hong Kong).

\noindent\hangindent 25pt I. Olkin, and H. Rubin, (1962). A
characterization of the Wishart distribution, \emph{Ann. Math.
Statist.}, 33, 1272-1280.

\noindent\hangindent 25pt J. Pauw, A. Bekker and J. J. J. Roux,
(2010). Densities of composite Weibullized generalized gamma
variables, \emph{S. Afr. Statist. J.}, 44, 17-42.

\noindent\hangindent 25pt G. O. Silva, E. M. M. Ortega, and G. M.
Cordeiro, (2010). The beta modified Weibull distribution,
\emph{Lifetime Data Anal.}, 16, 409-430.

\noindent\hangindent 25pt N. Singla, K. Jain, and S. K. Sharma,
(2012). The beta generalized Weibull distribution: Properties and applications, \emph{%
Reliability Eng. Sys. Safety}, 102, 5-15.

\noindent\hangindent 25pt S. M. Sunoj, and M. N. Linu, (2012).
Dynamic cumulative residual Renyi's entropy, \emph{Statistics},
46(1), 41-56.

\noindent\hangindent 25pt J. Van Niekerk, M. Arashi, A. Bekker and
J. J. J. Roux, (2013). Wishart kernel oriented generator
distributions with application, \emph{Technical Report},
University of Pretoria, South Africa.

\noindent \hangindent25pt K. Zagrafos and S. Nadarajah, (2005). Expression for R%
\'{e}nyi and Shannon entropies for multivariate distributions, \emph{%
Statist. Prob. Lett.}, 71, 71-84.

\end{document}